\documentclass[12pt]{article}
\usepackage{latexsym,amsfonts,amssymb}
\usepackage[dvips]{color}
\usepackage{colordvi,multicol}
\usepackage{amsmath}

\def \cal{\mathcal}

\newtheorem{thm}{Theorem}[section]
\newtheorem{cor}[thm]{Corollary}
\newtheorem{lem}[thm]{Lemma}
\newtheorem{pro}[thm]{Proposition}

\newtheorem{rem}[thm]{Remark}

\date{}
\begin{document}

\title{\bf The 4-D Gaussian Random Vector Maximum Conjecture and the 3-D Simplex Mean Width Conjecture}
 \author{Wei Sun$^{a,}\footnote{Corresponding author: wei.sun@concordia.ca}$,\ \ Ze-Chun Hu$^b$ and Guolie Lan$^c$\\ \\ \\
{\small $^a$ Department of Mathematics and Statistics, Concordia
University, Canada}\\ \\ {\small $^b$  College of Mathematics, Sichuan  University,  China }\\ \\
{\small $^c$ School of Economics and Statistics, Guangzhou University, China}\\ \\
}

\maketitle

\begin{abstract}
\vskip 0.3cm
\noindent We prove the four-dimensional Gaussian random vector maximum conjecture. This conjecture asserts that among all centered
Gaussian random vectors $X=(X_1,X_2,X_3,X_4)$ with $E[X_i^2]=1$, $1\le i\le 4$, the expectation $E[\max(X_1,X_2,X_3,X_4)]$ is maximal if and only if all  off-diagonal elements of the covariance matrix equal $-\frac{1}{3}$. As a direct consequence, we resolve the three-dimensional simplex mean width conjecture. This latter conjecture is a long-standing open problem in convex geometry, which asserts that among all simplices inscribed into the three-dimensional unit Euclidean ball the regular simplex has the maximal mean width.

\end{abstract}

\smallskip


\smallskip


\section{Introduction and main results}

The multivariate Gaussian distribution has played a fundamental role in probability and statistics (see \cite{T} and references therein). Although its study has a long history, there are remarkable new results which have been obtained in recent years. For example, Royen \cite{Royen} proved the Gaussian correlation inequality (cf. also Lata{\l}a and Matlak \cite{LM}), which has important applications in small ball probabilities (cf. Li \cite{Li99} and Shao \cite{Sh03}). As another example, based on stochastic calculus, Eldan \cite{E} gave a novel short proof of the celebrated Borell inequality and obtained an almost tight, two-sided, dimension-free robustness estimate
for the Gaussian noise stability
deficit. However, some unsolved problems still remain.   A major one is the Gaussian random vector maximum (GRVM) conjecture. This conjecture asserts that among all centered
Gaussian random vectors $X=(X_1,\dots,X_n)$ with $E[X_i^2]=1$, $1\le i\le n$, the expectation $E[\max(X_1,\dots,X_n)]$ is maximal if and only if all  off-diagonal elements of the covariance matrix equal $-\frac{1}{n-1}$, where $n\ge 4$. We would like to point out that the GRVM conjecture is much more difficult than another problem of order statistics of Gaussian random variables. If we add absolute values to $X_i$'s, then it is known that the expectation $E[\max(|X_1|,\dots,|X_n|)]$ attains maximum if and only $X_1,\dots,X_n$ are independent (see \v{S}id\'ak \cite{Si} and Gluskin \cite{Gl}).

Besides its own interest in probability and statistics, the GRVM conjecture has a deep connection with convex geometry. It is known that the $n$-dimensional GRVM conjecture is equivalent to the $(n-1)$-dimensional simplex mean width (SMW) conjecture (see Litvak \cite{Litvak}). Given a convex body $K$ in $\mathbb{R}^{n}$, its support function and mean width are defined by
$$
h_K(u)=\max_{x\in K}\langle u,x\rangle,\ \ \ \ \ \ \ \ w(K)=2\int_{s^{n-1}}h_K(u)dm(u).
$$
Hereafter $\langle\cdot,\cdot\rangle$ and $\|\cdot\|$ denote the Euclidean inner product and norm, respectively, and $m$ is the normalized Lebesgue measure on the sphere $S^{n-1}$. The SMW conjecture asserts that among all simplices inscribed into the $n$-dimensional unit ball the regular simplex has the maximal mean width, where $n\ge 3$.

The SMW conjecture is a long-standing open problem in convex geometry. As pointed out by Gritzmann and Klee \cite[Section 9.10.2]{KL}, several authors (\cite {Gi52,Ba63,Ba65,We68}) have assumed the existence of a proof for the conjecture. However, we are not aware of any such proof. As mentioned in recent works by  Hug and Schneider \cite{HS} and by B\"or\"ozky and Schneider \cite{BS}, the $n$-dimensional SMW conjecture remains open for $n\ge 3$. In the information theory community, it was a general belief that the SMW  conjecture is true (cf. \cite{Ba63,Ba65,D,We}) and  the SMW  conjecture is of interest because of its connection with the famous simplex code conjecture (see \cite[Section 9.10.2]{KL} or \cite{Litvak} for an introduction). For more history about the GRVM conjecture and the SMW conjecture, we refer the reader to the beautiful survey paper of Litvak \cite{Litvak}.

In this paper, we will prove the four-dimensional GRVM conjecture and, equivalently, the three-dimensional SMW conjecture.  Let $(\Omega,{\cal F},P)$ be a probability space such that, for any nonnegative-definite symmetric real matrix $\Lambda$, there exists a Gaussian random vector $X$ on $(\Omega,{\cal F},P)$ whose covariance matrix is given by $\Lambda$. We denote by ${\cal G}$ the set of all four-dimensional centered
Gaussian random vectors $X=(X_1,X_2,X_3,X_4)$ with $E[X_1^2]=E[X_2^2]=E[X_3^2]=E[X_4^2]=1$. Hereafter $E$ denotes the expectation with respect to $P$. Define
$$
{\cal M}(X)=E[\max(X_1,X_2,X_3,X_4)], \ \ \ \ X\in {\cal G}.
$$
For $X\in{\cal G}$, we denote its covariance matrix by $\Lambda^X$ or simply by $\Lambda$ if no confusion arises. We use $\bar{\cal S}$ and ${\cal S}$ to denote respectively the set of all nonnegative-definite
and positive-definite $4\times 4$ symmetric real
matrices with all diagonal elements equal to 1. For $\Lambda\in{\bar{\cal S}}$,
take an $X^{\Lambda}\in{\cal G}$ with covariance matrix $\Lambda$ and define
$$
F(\Lambda)={\cal M}(X^{\Lambda}).
$$

Now we state the main results of this paper.
\begin{thm}\label{main} The function ${\cal M}$ defined on ${\cal G}$ attains its maximum at $X$ if and only if
$\Lambda^X_{kl}=-\frac{1}{3}$ for all $k<l$.
\end{thm}

As a direct consequence of Theorem \ref{main}, we conclude that the three-dimensional SMW conjecture is true.

\begin{thm}
Among all simplices inscribed into the three-dimensional unit Euclidean ball the regular simplex has the maximal mean width.
\end{thm}

The SMW conjecture is a  purely geometric problem. We resolve the three-dimensional SMW conjecture by combining probability theory and convex geometry. First, we use a probabilistic method to derive the first order partial derivative formula for the function $F(\Lambda)$. To the best of our knowledge, this formula is unknown before. Based on this novel formula, we analyze the maximum point of the function ${\cal M}(X)$. By applying geometric techniques, we discover an intrinsic relation satisfied by the maximum point (see (\ref{WEI}) below). This important relation is described in terms of the dihedral angles of the tetrahedron $\{X_1,X_2,X_3,X_4\}$. We will show that the regular tetrahedron is the unique tetrahedron that satisfies the relation. To this end, we introduce a special function $H$ (see (\ref{WW5}) below) and reduce the three-dimensional SMW conjecture to the monotonicity of  $H$. We will use purely analytic techniques to prove that $H$ is a strictly decreasing function. This result has independent interest and some exquisitely delicate inequalities are involved in the proof. We hope the methods developed in this paper can be further developed so as to investigate the higher dimensional GRVM and SMW conjectures.

The remainder of this paper is organized as follows. In Section 2, we derive the first order partial derivative formula for $F(\Lambda)$. By virtue of this formula, we derive in Section 3 the second order partial derivative formula for $F(\Lambda)$ and the explicit formula for $F(\Lambda)$ itself. In Section 4, we prove Theorem \ref{main}. The proof relies heavily on  the monotonicity of the special function $H$ mentioned above. In Section 5, we prove  that $H$ is a strictly decreasing function. In Section 6, we make some remarks including non-concavity of the function $F(\Lambda)$ and the lower and upper bounds for
the function ${\cal M}(X)$. Some auxiliary results are given in the Appendix, Section 7.

\section{First order partial derivatives of $F(\Lambda)$}

Under the embedding map $I(\Lambda):=(\Lambda_{12}, \Lambda_{13},\Lambda_{14},\Lambda_{23},\Lambda_{24},\Lambda_{34})$, $\bar{\cal S}$ and ${\cal S}$ become respectively closed and open  convex subsets of $\mathbb{R}^{6}$. Let $\Lambda\in{\cal S}$. We define
\begin{equation}\label{ffff}
f_{\Lambda}(x)=\frac{1}{\sqrt{(2\pi)^4{\rm det}(\Lambda)}}e^{-\frac{1}{2}x^T{\Lambda}^{-1}x},\ \ x\in \mathbb{R}^4.
\end{equation}
Hereafter we use $^T$ to denote the transpose of a vector or a matrix.

\begin{lem}\label{Gauss2} For $\Lambda\in{\cal S}$, we have
\begin{eqnarray}\label{F12}
\frac{\partial F}{\partial\Lambda_{12}}(\Lambda)=-\int_{-\infty}^{\infty}\int_{-\infty}^t\int_{-\infty}^tf_{\Lambda}(t,t,x_3,x_4)dx_3dx_4dt.
\end{eqnarray}
Similar formulas hold for all $k<l$. Hence
$$
\frac{\partial F}{\partial\Lambda_{kl}}(\Lambda)<0,\ \ \ \ \forall k<l.
$$
\end{lem}

\noindent {\bf Proof.}\ \ Let $X=(X_1,X_2,X_3,X_4)\in {\cal G}$. Note that
$$
\lim_{t\uparrow\infty}tP\{\max(X_1,X_2,X_3,X_4)>t\}\le \sum_{k=1}^4\lim_{t\uparrow\infty}tP\{X_k>t\}=0.
$$
Then, we have
\begin{eqnarray}\label{W0}
&&E[\max(X_1,X_2,X_3,X_4)]\nonumber\\
&=&\int_{-\infty}^{\infty}tdP\{\max(X_1,X_2,X_3,X_4)<t\}\nonumber\\
&=&-\int_{0}^{\infty}tdP\{\max(X_1,X_2,X_3,X_4)>t\}+\int_{-\infty}^{0}tdP\{\max(X_1,X_2,X_3,X_4)<t\}\nonumber\\
&=&\int_{0}^{\infty}P\{\max(X_1,X_2,X_3,X_4)>t\}dt-\int_{-\infty}^{0}P\{\max(X_1,X_2,X_3,X_4)<t\}dt\nonumber\\
&=&\int_{0}^{\infty}\left[1-P\{\max(X_1,X_2,X_3,X_4)<t\}\right]dt-\int_0^{\infty}P\{\max(X_1,X_2,X_3,X_4)<-t\}dt\nonumber\\
&=&\int_{0}^{\infty}\left[1-P\{\max(X_1,X_2,X_3,X_4)<t\}-P\{\min(X_1,X_2,X_3,X_4)>t\}\right]dt.
\end{eqnarray}

We have
\begin{eqnarray*}
f_{\Lambda}(x)=\frac{1}{(2\pi)^4}\int_{\mathbb{R}^4} e^{-is^Tx-\frac{1}{2}s^T{\Lambda}s}ds=
\frac{1}{(2\pi)^4}\int_{\mathbb{R}^4} e^{-is^Tx-\frac{1}{2}\sum_{k,l=1}^4{\Lambda}_{kl}s_ks_l}ds,\ \ x\in \mathbb{R}^4,
\end{eqnarray*}
and
\begin{eqnarray}\label{Fou}
\frac{\partial f_{\Lambda}}{\partial{\Lambda}_{kl}}=\frac{\partial^2 f_{\Lambda}}{\partial x_k\partial x_l},\ \ \ \ \forall k<l.
\end{eqnarray}
By (\ref{W0}), we get
{\small\begin{eqnarray*}
F(\Lambda)&=&\int_{0}^{\infty}\left[1-P\{\max(X^{\Lambda}_1,\dots, X^{\Lambda}_4)<t\}-P\{\min(X^{\Lambda}_1,\dots, X^{\Lambda}_4)>t\}\right]dt\\
&=&\int_0^{\infty}\left[1-\int_{-\infty}^t\cdots \int_{-\infty}^tf_{\Lambda}(x_1,\dots,x_4)
dx_1\cdots dx_4-\int_t^{\infty}\cdots \int_t^{\infty}f_{\Lambda}(x_1,\dots,x_4)dx_1\cdots dx_4\right]dt.
\end{eqnarray*}}
Then, we obtain by (\ref{ffff}) and (\ref{Fou}) that
\begin{eqnarray*}
&&\frac{\partial F}{\partial\Lambda_{12}}(\Lambda)\nonumber\\
&=&-\int_0^{\infty}\left[\int_{-\infty}^t\cdots \int_{-\infty}^t
\frac{\partial f_{\Lambda}}{\partial{\Lambda}_{12}}dx_1\cdots dx_4+\int_t^{\infty}\cdots \int_t^{\infty}\frac{\partial f_{\Lambda}}{\partial{\Lambda}_{12}}dx_1\cdots dx_4\right]dt\nonumber\\
&=&-\int_0^{\infty}\left[\int_{-\infty}^t\cdots \int_{-\infty}^t\frac{\partial^2 f_{\Lambda}}{\partial x_1\partial x_2}dx_1\cdots dx_4+\int_t^{\infty}\cdots \int_t^{\infty}\frac{\partial^2 f_{\Lambda}}{\partial x_1\partial x_2}dx_1\cdots dx_4\right]dt\nonumber\\
&=&-\int_0^{\infty}\left[\int_{-\infty}^t\int_{-\infty}^tf_{\Lambda}(t,t,x_3,x_4)dx_3dx_4+\int_t^{\infty} \int_t^{\infty}f_{\Lambda}(t,t,x_3,x_4)dx_3dx_4\right]dt\nonumber\\
&=&-\int_0^{\infty}\int_{-\infty}^t\int_{-\infty}^tf_{\Lambda}(t,t,x_3,x_4)dx_3dx_4dt+\int_{-\infty}^0\int_{-\infty}^t \int_{-\infty}^tf_{\Lambda}(-t,-t,-x_3,-x_4)dx_3dx_4dt\nonumber\\
&=&-\int_{-\infty}^{\infty}\int_{-\infty}^t\int_{-\infty}^tf_{\Lambda}(t,t,x_3,x_4)dx_3dx_4dt\nonumber\\
&<&0.
\end{eqnarray*}
\hfill\fbox

\begin{rem}
Lemma \ref{Gauss2} can apparently be extended to higher dimensional cases.
\end{rem}

Suppose $\Lambda\in{\bar{\cal S}}$. We choose an $X^{\Lambda}=(X^{\Lambda}_1,X^{\Lambda}_2,X^{\Lambda}_3,X^{\Lambda}_4)\in{\cal G}$ with covariance matrix $\Lambda$. Regard
$X^{\Lambda}_1,X^{\Lambda}_2,X^{\Lambda}_3,X^{\Lambda}_4$ as vectors of the space $L^2(\Omega,{\cal F},P)$ and denote by ${\cal V}$ the volume of the tetrahedron with vertexes $\{X^{\Lambda}_1,X^{\Lambda}_2,X^{\Lambda}_3,X^{\Lambda}_4\}$. Define
$$
{\tilde A}=6\sqrt{2}{\cal V}.
$$
Note that $({X^{\Lambda}_1-X^{\Lambda}_2}$, ${X^{\Lambda}_3-X^{\Lambda}_2}$, ${X^{\Lambda}_4-X^{\Lambda}_2})$ is a 3-dimensional centered
Gaussian random vector with covariance matrix
\begin{eqnarray}\label{extend}
\Sigma^{\Lambda}_2:=\left(
\begin{array}{ccc}
2(1-{\Lambda}_{12})&{1-{\Lambda}_{12}+{\Lambda}_{13}-{\Lambda}_{23}} & {1-{\Lambda}_{12}+{\Lambda}_{14}-{\Lambda}_{24}}\\
{1-{\Lambda}_{12}+{\Lambda}_{13}-{\Lambda}_{23}} &2(1-{\Lambda}_{23})&{1-{\Lambda}_{23}-{\Lambda}_{24}+{\Lambda}_{34}}\\
{1-{\Lambda}_{12}+{\Lambda}_{14}-{\Lambda}_{24}}&{1-{\Lambda}_{23}-{\Lambda}_{24}+{\Lambda}_{34}}&2(1-{\Lambda}_{24})
\end{array}
\right).
\end{eqnarray}
By the volume formula of a tetrahedron, we get
$$
{\cal V}=\frac{1}{6}\sqrt{{\rm det}(\Sigma^{\Lambda}_2)}.
$$
Hence
\begin{equation}\label{AAAAA}
{\tilde A}=\sqrt{2{\rm det}(\Sigma^{\Lambda}_2)}.
\end{equation}
Let $\Sigma^{\Lambda}_k$ be the covariance matrix of the Gaussian random vector $({X^{\Lambda}_{l_1}-X^{\Lambda}_k}$, ${X^{\Lambda}_{l_2}-X^{\Lambda}_k}$, ${X^{\Lambda}_{l_3}-X^{\Lambda}_k})$, where $l_1<l_2<l_3$ is the increasing arrangement of  $\{1,2,3,4\}\setminus\{k\}$. Similar to (\ref{AAAAA}), we can show that
$$
{\tilde A}=\sqrt{2{\rm det}(\Sigma^{\Lambda}_k)},\ \ \ \ k\in\{1,2,3,4\}.
$$

Define
$$
(\Lambda')_{kl}=1-{\Lambda}_{kl},\ \ \ \ k<l,
$$
and
\begin{eqnarray}\label{qqq}
(\tilde{\Lambda})_{12}&=&(1-{\Lambda}_{12}+{\Lambda}_{14}-{\Lambda}_{24})(1-{\Lambda}_{12}+{\Lambda}_{13}-{\Lambda}_{23})
-2(1-{\Lambda}_{12})(1-{\Lambda}_{23}-{\Lambda}_{24}+{\Lambda}_{34})\nonumber\\
&=&(\Lambda')^2_{12}-(\Lambda')_{12}[(\Lambda')_{13}+(\Lambda')_{14}+(\Lambda')_{23}
+(\Lambda')_{24}-2(\Lambda')_{34}]\nonumber\\
&&+[(\Lambda')_{13}-(\Lambda')_{23}]\cdot[(\Lambda')_{14}-(\Lambda')_{24}].
\end{eqnarray}
Similarly, we define $(\tilde{\Lambda})_{kl}$ for $k<l$ by replacing $1$ in (\ref{qqq}) with $k$, $2$ in (\ref{qqq}) with $l$,
$3$ in (\ref{qqq}) with $\min\{s:1\le s\le 4, s\not=k,l\}$, and 4 in (\ref{qqq}) with
$\max\{s:1\le s\le 4, s\not=k,l\}$.

Suppose $\Lambda\in {\cal S}$. We have
\begin{eqnarray}\label{hao0}
&&\left(\begin{array}{cccc}
1&0&0&0\\
1&1&1&1\\
0&0&1&0\\
0&0&0&1
\end{array}\right)\left(\begin{array}{cccc}
(\Lambda^{-1})_{11}&(\Lambda^{-1})_{12}&(\Lambda^{-1})_{13}&(\Lambda^{-1})_{14}\\
(\Lambda^{-1})_{21}&(\Lambda^{-1})_{22}&(\Lambda^{-1})_{23}&(\Lambda^{-1})_{24}\\
(\Lambda^{-1})_{31}&(\Lambda^{-1})_{32}&(\Lambda^{-1})_{33}&(\Lambda^{-1})_{34}\\
(\Lambda^{-1})_{41}&(\Lambda^{-1})_{42}&(\Lambda^{-1})_{43}&(\Lambda^{-1})_{44}\nonumber
\end{array}\right)\left(\begin{array}{cccc}
1&1&0&0\\
0&1&0&0\\
0&1&1&0\\
0&1&0&1
\end{array}\right)\\
&=&\left(\begin{array}{cccc}
(\Lambda^{-1})_{11}&\sum_{l=1}^4(\Lambda^{-1})_{1l}&(\Lambda^{-1})_{13}&(\Lambda^{-1})_{14}\\
\sum_{l=1}^4(\Lambda^{-1})_{1l}&\sum_{k,l=1}^4(\Lambda^{-1})_{kl}&\sum_{l=1}^4(\Lambda^{-1})_{3l}&\sum_{l=1}^4(\Lambda^{-1})_{4l}\\
(\Lambda^{-1})_{31}&\sum_{l=1}^4(\Lambda^{-1})_{3l}&(\Lambda^{-1})_{33}&(\Lambda^{-1})_{34}\\
(\Lambda^{-1})_{41}&\sum_{l=1}^4(\Lambda^{-1})_{4l}&(\Lambda^{-1})_{43}&(\Lambda^{-1})_{44}
\end{array}\right).
\end{eqnarray}
Then,
$$
\sum_{k,l=1}^4(\Lambda^{-1})_{kl}>0.
$$
Define
$$
A=\sqrt{\sum_{k,l=1}^4(\Lambda^{-1})_{kl}}.
$$
By (\ref{hao0}), we get
{\small\begin{eqnarray}\label{hao}
&&\left(\begin{array}{cccc}
(\Lambda^{-1})_{11}&\sum_{l=1}^4(\Lambda^{-1})_{1l}&(\Lambda^{-1})_{13}&(\Lambda^{-1})_{14}\nonumber\\
\sum_{l=1}^4(\Lambda^{-1})_{1l}&A^2&\sum_{l=1}^4(\Lambda^{-1})_{3l}&\sum_{l=1}^4(\Lambda^{-1})_{4l}\nonumber\\
(\Lambda^{-1})_{31}&\sum_{l=1}^4(\Lambda^{-1})_{3l}&(\Lambda^{-1})_{33}&(\Lambda^{-1})_{34}\nonumber\\
(\Lambda^{-1})_{41}&\sum_{l=1}^4(\Lambda^{-1})_{4l}&(\Lambda^{-1})_{43}&(\Lambda^{-1})_{44}\nonumber\\
\end{array}\right)^{-1}\nonumber\\
&=&\left(\begin{array}{cccc}
1&1&0&0\\
0&1&0&0\\
0&1&1&0\\
0&1&0&1
\end{array}\right)^{-1}\Lambda\left(\begin{array}{cccc}
1&0&0&0\nonumber\\
1&1&1&1\nonumber\\
0&0&1&0\nonumber\\
0&0&0&1
\end{array}\right)^{-1}\\
&=&\left(\begin{array}{cccc}
1&-1&0&0\nonumber\\
0&1&0&0\nonumber\\
0&-1&1&0\nonumber\\
0&-1&0&1
\end{array}\right)\Lambda\left(\begin{array}{cccc}
1&0&0&0\nonumber\\
-1&1&-1&-1\nonumber\\
0&0&1&0\nonumber\\
0&0&0&1
\end{array}\right)\nonumber\\
&=&\left(\begin{array}{cccc}
2(1-{\Lambda}_{12})&{\Lambda}_{12}-1&1-{\Lambda}_{12}+{\Lambda}_{13}-{\Lambda}_{23}
&1-{\Lambda}_{12}+{\Lambda}_{14}-{\Lambda}_{24}\\
{\Lambda}_{12}-1&1&{\Lambda}_{23}-1&{\Lambda}_{24}-1\\
1-{\Lambda}_{12}+{\Lambda}_{13}-{\Lambda}_{23}&{\Lambda}_{23}-1&2(1-{\Lambda}_{23})
&1-{\Lambda}_{23}-{\Lambda}_{24}+{\Lambda}_{34}\\
1-{\Lambda}_{12}+{\Lambda}_{14}-{\Lambda}_{24}&{\Lambda}_{24}-1&1-{\Lambda}_{23}-{\Lambda}_{24}+{\Lambda}_{34}
&2(1-{\Lambda}_{24})\end{array}\right).
\end{eqnarray}}
Then, we obtain by (\ref{AAAAA}) and (\ref{hao}) that
{\begin{eqnarray}\label{dis3}
A^2=\frac{{\rm det}(\Sigma^{\Lambda}_2)}{{\rm det}({\Lambda})}=\frac{{\tilde A}^2}{2{\rm det}({\Lambda})}.
\end{eqnarray}

\begin{thm}\label{Gauss33} For $\Lambda\in{\cal S}$, we have
\begin{equation}\label{June0}\frac{\partial F}{\partial\Lambda_{kl}}(\Lambda)
=-\frac{1}{4\sqrt{\pi^3(\Lambda')_{kl}}}\arccos\left(\frac{(\tilde{\Lambda})_{kl}}{\sqrt{({\Lambda}')_{kl}{\tilde A}^2+(\tilde{\Lambda})^2_{kl}}}\right)
,\ \ \ \ \forall k<l.
\end{equation}
\end{thm}

\noindent {\bf Proof.}\ \ Without loss of generality, we prove (\ref{June0}) for the case that $k=1,l=2$. By (\ref{F12}), we get
{\small\begin{eqnarray}\label{nan}
&&-\frac{\partial F}{\partial\Lambda_{12}}(\Lambda)\nonumber\\
&=&\int_{-\infty}^{\infty}\int_{-\infty}^t\int_{-\infty}^tf_{\Lambda}(t,t,x_3,x_4)dx_3dx_4dt\nonumber\\
&=&\int_{-\infty}^{\infty}\int_{-\infty}^0\int_{-\infty}^0f_{\Lambda}(t,t,x_3+t,x_4+t)dx_3dx_4dt\nonumber\\
&=&\int_{-\infty}^0\int_{-\infty}^0\int_{-\infty}^{\infty}f_{\Lambda}(t,t,x_3+t,x_4+t)dt dx_3dx_4\nonumber\\
&=&\frac{1}{\sqrt{(2\pi)^{4}{\rm det}({\Lambda})}}\int_{-\infty}^0\int_{-\infty}^0\int_{-\infty}^{\infty}\nonumber\\
&&\ \ \ \ e^{-\frac{1}{2}\left\{\sum_{k,l=1}^4(\Lambda^{-1})_{kl}t^2+2\left[\sum_{l=1}^4(\Lambda^{-1})_{3l}x_3+\sum_{l=1}^4(\Lambda^{-1})_{4l}x_4\right]t
+\left[(\Lambda^{-1})_{33}x^2_3+2(\Lambda^{-1})_{34}x_3x_4+(\Lambda^{-1})_{44}x^2_4\right]\right\}}dt dx_3 dx_4\nonumber\\
&=&\frac{1}{\sqrt{(2\pi)^{4}{\rm det}({\Lambda})}}\int_{-\infty}^0\int_{-\infty}^0\int_{-\infty}^{\infty}
 e^{-\frac{1}{2}\left\{\sqrt{\sum_{k,l=1}^4(\Lambda^{-1})_{kl}}t+\frac{\sum_{l=1}^4(\Lambda^{-1})_{3l}x_3+
 \sum_{l=1}^4(\Lambda^{-1})_{4l}x_4}{\sqrt{\sum_{k,l=1}^4(\Lambda^{-1})_{kl}}}\right\}^2}dt\nonumber\\
&&\ \ \ \ \ \ \ \ \ \ \ \ \ \ \ \ \ \ \ \ e^{-\frac{1}{2}\left\{\left[(\Lambda^{-1})_{33}x^2_3
+2(\Lambda^{-1})_{34}x_3x_4+(\Lambda^{-1})_{44}x^2_4\right]-\frac{\left[\sum_{l=1}^4(\Lambda^{-1})_{3l}x_3
+\sum_{l=1}^4(\Lambda^{-1})_{4l}x_4\right]^2}{{\sum_{k,l=1}^4(\Lambda^{-1})_{kl}}}\right\}}dx_3 dx_4\nonumber\\
&=&\frac{1}{\sqrt{(2\pi)^{3}{\rm det}({\Lambda})\sum_{k,l=1}^4(\Lambda^{-1})_{kl}}}\int_{-\infty}^0\int_{-\infty}^0\nonumber\\
&&\ \ \ \ \ \ \ \ \ \ \ \ \ \ \ \ \ \ \ \ e^{-\frac{1}{2}\left\{\left[(\Lambda^{-1})_{33}x^2_3
+2(\Lambda^{-1})_{34}x_3x_4+(\Lambda^{-1})_{44}x^2_4\right]-\frac{\left[\sum_{l=1}^4(\Lambda^{-1})_{3l}x_3
+\sum_{l=1}^4(\Lambda^{-1})_{4l}x_4\right]^2}{{\sum_{k,l=1}^4(\Lambda^{-1})_{kl}}}\right\}}dx_3 dx_4\nonumber\\
&=&\frac{1}{\sqrt{(2\pi)^{3}{\rm det}({\Lambda})\sum_{k,l=1}^4(\Lambda^{-1})_{kl}}}\int_0^{\infty}\int_0^{\infty}\nonumber\\
&&\ \ \ \ \ \ \ \ \ \ \ \ \ \ \ \ \ \ \ \ e^{-\frac{1}{2}\left\{\left[(\Lambda^{-1})_{33}x^2_3
+2(\Lambda^{-1})_{34}x_3x_4+(\Lambda^{-1})_{44}x^2_4\right]-\frac{\left[\sum_{l=1}^4(\Lambda^{-1})_{3l}x_3
+\sum_{l=1}^4(\Lambda^{-1})_{4l}x_4\right]^2}{{\sum_{k,l=1}^4(\Lambda^{-1})_{kl}}}\right\}}dx_3 dx_4\nonumber\\
&:=&\frac{1}{\sqrt{(2\pi)^{3}{\rm det}({\Lambda})\sum_{k,l=1}^4(\Lambda^{-1})_{kl}}}\int_0^{\infty}\int_0^{\infty}
e^{-\frac{a_1x_3^2+b_1x_4^2+2c_1x_3x_4}{2A^2}}dx_3dx_4,
\end{eqnarray}}
where
$$
a_1:=(\Lambda^{-1})_{33}A^2-\left(\sum_{l=1}^4(\Lambda^{-1})_{3l}\right)^2,\ \ \ \ b_1:=(\Lambda^{-1})_{44}A^2-\left(\sum_{l=1}^4(\Lambda^{-1})_{4l}\right)^2,
$$
and
\begin{equation}\label{215}
c_1:=(\Lambda^{-1})_{34}A^2-\left(\sum_{l=1}^4(\Lambda^{-1})_{3l}\right)\left(\sum_{l=1}^4(\Lambda^{-1})_{4l}\right).
\end{equation}

Set
$$
t^*=-\frac{\sum_{l=1}^4(\Lambda^{-1})_{3l}x_3+
 \sum_{l=1}^4(\Lambda^{-1})_{4l}x_4}{A^2}.
$$
Then,
\begin{eqnarray*}
&&\frac{a_1x_3^2+b_1x_4^2+2c_1x_3x_4}{A^2}\\
&=&[(\Lambda^{-1})_{33}x^2_3
+2(\Lambda^{-1})_{34}x_3x_4+(\Lambda^{-1})_{44}x^2_4]-\frac{[\sum_{l=1}^4(\Lambda^{-1})_{3l}x_3
+\sum_{l=1}^4(\Lambda^{-1})_{4l}x_4]^2}{{\sum_{k,l=1}^4(\Lambda^{-1})_{kl}}}\\
&=&(t^*,t^*,x_3+t^*,x_4+t^*)(\Lambda^{-1})(t^*,t^*,x_3+t^*,x_4+t^*)^T.
\end{eqnarray*}
Note that $(x_3,x_4)\not=0$ implies  $(t^*,t^*,x_3+t^*,x_4+t^*)\not=0$. Hence, we have
\begin{equation}\label{a1b1}
a_1b_1-c_1^2>0.
\end{equation}
By (\ref{nan}) and Lemma \ref{lem71} of the Appendix, we get

{\small\begin{eqnarray}\label{OO}
\int_{-\infty}^{\infty}\int_{-\infty}^t\int_{-\infty}^tf_{\Lambda}(t,t,x_3,x_4)dx_3dx_4dt
=\frac{A}{\sqrt{(2\pi)^3{\rm det}({\Lambda})(a_1b_1-c^2_1)}}\arccos\left(\frac{c_1}{\sqrt{a_1b_1}}\right).
\end{eqnarray}}

Define
\begin{equation}\label{IIII}
\Pi=\left(\begin{array}{ccc}
A^2&\sum_{l=1}^4(\Lambda^{-1})_{3l}&\sum_{l=1}^4(\Lambda^{-1})_{4l}\\
\sum_{l=1}^4(\Lambda^{-1})_{3l}&(\Lambda^{-1})_{33}&(\Lambda^{-1})_{34}\\
\sum_{l=1}^4(\Lambda^{-1})_{4l}&(\Lambda^{-1})_{34}&(\Lambda^{-1})_{44}
\end{array}\right).
\end{equation}
Then, we obtain by (\ref{hao0}) and (\ref{hao}) that
\begin{equation}\label{17}
{\rm det}({\Lambda})\cdot{\rm det}(\Pi)=\frac{{\rm det}(\Pi)}{{\rm det}(\Lambda^{-1})}=2(1-{\Lambda}_{12}).
\end{equation}
Thus,
{\small\begin{eqnarray}\label{dis}
&&{\rm det}({\Lambda})(a_1b_1-c^2_1)\nonumber\\
&=&{\rm det}({\Lambda})\left[((\Lambda^{-1})_{33}(\Lambda^{-1})_{44}-(\Lambda^{-1})^2_{34})A^4\right.\nonumber\\
&&\left.-\left((\Lambda^{-1})_{33}
\left(\sum_{l=1}^4(\Lambda^{-1})_{4l}\right)^2 + (\Lambda^{-1})_{44}
\left(\sum_{l=1}^4(\Lambda^{-1})_{3l}\right)^2-2(\Lambda^{-1})_{34}
\left(\sum_{l=1}^4(\Lambda^{-1})_{3l}\right)\left(\sum_{l=1}^4(\Lambda^{-1})_{4l}\right)  \right)A^2\right]\nonumber\\
&=&A^2{\rm det}({\Lambda})\cdot{\rm det (\Pi)}\nonumber\\
&=&2(1-{\Lambda}_{12})A^2,
\end{eqnarray}}
which implies that
\begin{eqnarray}\label{BB}
\frac{A}{\sqrt{(2\pi)^3{\rm det}({\Lambda})(a_1b_1-c^2_1)}}
=\frac{1}{4\sqrt{\pi^3(1-{\Lambda}_{12})}}=\frac{1}{4\sqrt{\pi^3(\Lambda')_{12}}}.
\end{eqnarray}

Recall the block matrix inversion formula
$$
\left(\begin{array}{cc}
S&T\\
U&V
\end{array}\right)^{-1}=\left(\begin{array}{cc}
S^{-1}+S^{-1}T(V-US^{-1}T)^{-1}US^{-1}&-S^{-1}T(V-US^{-1}T)^{-1}\\
-(V-US^{-1}T)^{-1}US^{-1}&(V-US^{-1}T)^{-1}
\end{array}\right).
$$
By (\ref{hao}), we get
\begin{eqnarray}\label{18}
\Pi&=&\left[\left(\begin{array}{ccc}
1&{\Lambda}_{23}-1&{\Lambda}_{24}-1\nonumber\\
{\Lambda}_{23}-1 &2(1-{\Lambda}_{23})
&1-{\Lambda}_{23}-{\Lambda}_{24}+{\Lambda}_{34}\nonumber\\
{\Lambda}_{24}-1&1-{\Lambda}_{23}-{\Lambda}_{24}+{\Lambda}_{34}
&2(1-{\Lambda}_{24})\end{array}\right)\right.\\
&&\left. -\left(\begin{array}{c}
{\Lambda}_{12}-1\nonumber\\
1-{\Lambda}_{12}+{\Lambda}_{13}-{\Lambda}_{23}\nonumber\\
1-{\Lambda}_{12}+{\Lambda}_{14}-{\Lambda}_{24}\end{array}\right)[2(1-{\Lambda}_{12})]^{-1}\left(\begin{array}{c}
{\Lambda}_{12}-1\nonumber\\
1-{\Lambda}_{12}+{\Lambda}_{13}-{\Lambda}_{23}\nonumber\\
1-{\Lambda}_{12}+{\Lambda}_{14}-{\Lambda}_{24}\end{array}\right)^T\right]^{-1}.\nonumber\\
&&
\end{eqnarray}
Then, we obtain by (\ref{qqq}), (\ref{215}), (\ref{IIII}), (\ref{17}) and (\ref{18}) that
\begin{eqnarray}\label{dis2}
c_1&=&-{\rm det}(\Pi)\cdot(\Pi^{-1})_{32}\nonumber\\
&=&-\frac{2(1-{\Lambda}_{12})}{{\rm det}({\Lambda})}\cdot\left[(1-{\Lambda}_{23}-{\Lambda}_{24}+{\Lambda}_{34})\right.\nonumber\\
&&\ \ \ \ \left.-\frac{(1-{\Lambda}_{12}+{\Lambda}_{14}-{\Lambda}_{24})(1-{\Lambda}_{12}+{\Lambda}_{13}-{\Lambda}_{23})}
{2(1-{\Lambda}_{12})}\right]\nonumber\\
&=&\frac{1}
{{\rm det}({\Lambda})}[(1-{\Lambda}_{12}+{\Lambda}_{14}-{\Lambda}_{24})(1-{\Lambda}_{12}+{\Lambda}_{13}-{\Lambda}_{23})\nonumber\\
&&\ \ \ \ -2(1-{\Lambda}_{12})(1-{\Lambda}_{23}-{\Lambda}_{24}+{\Lambda}_{34})]\nonumber\\
&=&\frac{(\tilde{\Lambda})_{12}}
{{\rm det}({\Lambda})}.
\end{eqnarray}
Further, by (\ref{dis3}), (\ref{a1b1}), (\ref{dis}) and (\ref{dis2}), we get
\begin{eqnarray}\label{AA}
0&<&a_1b_1\nonumber\\
&=&(a_1b_1-c_1^2)+c_1^2\nonumber\\
&=&\frac{2{\rm det}({\Lambda}) ({\Lambda}')_{12}A^2+(\tilde{\Lambda})^2_{12}}{[{\rm det}({\Lambda})]^2}\nonumber\\
&=&\frac{({\Lambda}')_{12}{\tilde A}^2+(\tilde{\Lambda})^2_{12}}{[{\rm det}({\Lambda})]^2},
\end{eqnarray}
which together with (\ref{dis2}) implies that
\begin{eqnarray}\label{CC}
\arccos\left(\frac{c_1}{\sqrt{a_1b_1}}\right)=\arccos\left(\frac{(\tilde{\Lambda})_{12}}{\sqrt{({\Lambda}')_{12}{\tilde A}^2+(\tilde{\Lambda})^2_{12}}}\right).
\end{eqnarray}
Therefore, we obtain by (\ref{OO}), (\ref{BB}) and (\ref{CC}) that
\begin{eqnarray}\label{secondd}
&&\int_{-\infty}^{\infty}\int_{-\infty}^t\int_{-\infty}^tf_{\Lambda}(t,t,x_3,x_4)dx_3dx_4dt=\frac{1}{4\sqrt{\pi^3(\Lambda')_{12}}}
\arccos\left(\frac{(\tilde{\Lambda})_{12}}{\sqrt{({\Lambda}')_{12}{\tilde A}^2+(\tilde{\Lambda})^2_{12}}}\right).\ \ \ \ \ \ \ \
\end{eqnarray}
\hfill\fbox

\section{Second order partial derivatives and formula of $F(\Lambda)$}\setcounter{equation}{0}
\begin{thm}\label{Gauss44} For $\Lambda\in{\cal S}$, we have
{\small\begin{eqnarray}\label{61}
&&\frac{\partial^2 F}{\partial\Lambda_{12}\partial\Lambda_{12}}(\Lambda)\nonumber\\
&=&-\frac{1}{8\sqrt{\pi^3(\Lambda')^3_{12}}}
\arccos\left(\frac{(\tilde{\Lambda})_{12}}{\sqrt{({\Lambda}')_{12}{\tilde A}^2+(\tilde{\Lambda})^2_{12}}}\right)
+\frac{1}{4\sqrt{\pi^3}(\Lambda')_{12}{\tilde A}}\nonumber\\
&&\ \ \cdot\left[\frac{(\tilde{\Lambda})_{13}[(\Lambda')_{12}+(\Lambda')_{23}-(\Lambda')_{13}]}{4(\Lambda')_{12}(\Lambda')_{23}
-[(\Lambda')_{12}+(\Lambda')_{23}-(\Lambda')_{13}]^2}
+\frac{(\tilde{\Lambda})_{14}[(\Lambda')_{12}+(\Lambda')_{24}-(\Lambda')_{14}]}{4(\Lambda')_{12}(\Lambda')_{24}
-[(\Lambda')_{12}+(\Lambda')_{24}-(\Lambda')_{14}]^2}\right],\nonumber\\
&&
\end{eqnarray}}
$$\frac{\partial^2 F}{\partial\Lambda_{12}\partial\Lambda_{13}}(\Lambda)\\
=-\frac{(\tilde{\Lambda})_{23}}{2\sqrt{\pi^3}{\tilde A}\{4(\Lambda')_{12}(\Lambda')_{23}
-[(\Lambda')_{12}+(\Lambda')_{23}-(\Lambda')_{13}]^2\}},
$$
and
$$\frac{\partial^2 F}{\partial\Lambda_{12}\partial\Lambda_{34}}(\Lambda)\\
=-\frac{1}{2\sqrt{\pi^3}{\tilde A}}.
$$
Similar formulas hold for all $k<l$ if we replace $1$ with $k$, $2$  with $l$,
$3$  with $\min\{s:1\le s\le 4, s\not=k,l\}$, and 4  with
$\max\{s:1\le s\le 4, s\not=k,l\}$.
\end{thm}

\noindent {\bf Proof.}\ \ By (\ref{F12}), (\ref{extend})--(\ref{qqq}), (\ref{secondd}) and (\ref{symm}) of the Appendix, we get
{\small\begin{eqnarray*}
&&-\frac{\partial^2 F}{\partial\Lambda_{12}\partial\Lambda_{12}}(\Lambda)\nonumber\\
&=&\frac{\partial}{\partial{\Lambda}_{12}}\left\{\int_{-\infty}^{\infty}\int_{-\infty}^t\int_{-\infty}^tf_{\Lambda}(t,t,x_3,x_4)dx_3dx_4dt\right\}\nonumber\\
&=&\frac{1}{8\sqrt{\pi^3(\Lambda')^3_{12}}}
\arccos\left(\frac{(\tilde{\Lambda})_{12}}{\sqrt{({\Lambda}')_{12}{\tilde A}^2+(\tilde{\Lambda})^2_{12}}}\right)\nonumber\\
&&-\frac{\sqrt{({\Lambda}')_{12}{\tilde A}^2
+(\tilde{\Lambda})^2_{12}}}{4\sqrt{\pi^3}(\Lambda')_{12}{\tilde A}}
\frac{\partial}{\partial{\Lambda}_{12}}\left(\frac{(\tilde{\Lambda})_{12}}
{\sqrt{({\Lambda}')_{12}{\tilde A}^2+(\tilde{\Lambda})^2_{12}}}\right)\nonumber\\
&=&\frac{1}{8\sqrt{\pi^3(\Lambda')^3_{12}}}
\arccos\left(\frac{(\tilde{\Lambda})_{12}}{\sqrt{({\Lambda}')_{12}{\tilde A}^2+(\tilde{\Lambda})^2_{12}}}\right)\nonumber\\
&&-\frac{1}{4\sqrt{\pi^3}(\Lambda')_{12}{\tilde A}}\cdot\{-2(\Lambda')_{12}+(\Lambda')_{13}+(\Lambda')_{14}
+(\Lambda')_{23}+(\Lambda')_{24}-2(\Lambda')_{34}\nonumber\\
&&\ \ \left.+\frac{(\tilde{\Lambda})_{12}[(\Lambda')_{13}+(\Lambda')_{23}-(\Lambda')_{12}]}{4(\Lambda')_{12}(\Lambda')_{23}
-[(\Lambda')_{12}+(\Lambda')_{23}-(\Lambda')_{13}]^2}
+\frac{(\tilde{\Lambda})_{12}[(\Lambda')_{14}+(\Lambda')_{24}-(\Lambda')_{12}]}{4(\Lambda')_{12}(\Lambda')_{24}
-[(\Lambda')_{12}+(\Lambda')_{24}-(\Lambda')_{14}]^2}\right\}\nonumber\\
&=&\frac{1}{8\sqrt{\pi^3(\Lambda')^3_{12}}}
\arccos\left(\frac{(\tilde{\Lambda})_{12}}{\sqrt{({\Lambda}')_{12}{\tilde A}^2+(\tilde{\Lambda})^2_{12}}}\right)
-\frac{1}{4\sqrt{\pi^3}(\Lambda')_{12}{\tilde A}}\nonumber\\
&&\ \ \cdot\left[\frac{(\tilde{\Lambda})_{13}[(\Lambda')_{12}+(\Lambda')_{23}-(\Lambda')_{13}]}{4(\Lambda')_{12}(\Lambda')_{23}
+[(\Lambda')_{12}+(\Lambda')_{23}-(\Lambda')_{13}]^2}
+\frac{(\tilde{\Lambda})_{14}[(\Lambda')_{12}+(\Lambda')_{24}-(\Lambda')_{14}]}{4(\Lambda')_{12}(\Lambda')_{24}
-[(\Lambda')_{12}+(\Lambda')_{24}-(\Lambda')_{14}]^2}\right].\nonumber\\
&&
\end{eqnarray*}}

By (\ref{IIII}), (\ref{17}) and (\ref{18}), we have that
\begin{eqnarray}\label{test1}
&&(\Lambda^{-1})_{44}A^2-\left(\sum_{l=1}^4(\Lambda^{-1})_{4l}\right)^2\nonumber\\
&=&{\rm det}(\Pi)\cdot(\Pi^{-1})_{22}\nonumber\\
&=&\frac{2(1-{\Lambda}_{12})}{{\rm det}({\Lambda})}\cdot\left[2(1-{\Lambda}_{23})-\frac{(1-{\Lambda}_{12}+{\Lambda}_{13}-{\Lambda}_{23})^2}
{2(1-{\Lambda}_{12})}\right]\nonumber\\
&=&\frac{1}
{{\rm det}({\Lambda})}[4(1-{\Lambda}_{12})(1-{\Lambda}_{23})-(1-{\Lambda}_{12}+{\Lambda}_{13}-{\Lambda}_{23})^2]\nonumber\\
&=&\frac{1}
{{\rm det}({\Lambda})}\left[4(\Lambda')_{12}(\Lambda')_{23}-[(\Lambda')_{12}+(\Lambda')_{23}-(\Lambda')_{13}]^2
\right].
\end{eqnarray}
By (\ref{dis2}) and symmetry, we get
\begin{eqnarray}\label{test2}
(\Lambda^{-1})_{14}A^2-\left(\sum_{l=1}^4(\Lambda^{-1})_{1l}\right)\left(\sum_{l=1}^4(\Lambda^{-1})_{4l}\right)=\frac{(\tilde{\Lambda})_{23}}
{{\rm det}({\Lambda})}.
\end{eqnarray}
Then, we obtain by (\ref{ffff}), (\ref{F12}), (\ref{Fou}), (\ref{dis3}), (\ref{test1}) and (\ref{test2}) that
\begin{eqnarray*}
&&-\frac{\partial^2 F}{\partial\Lambda_{12}\partial\Lambda_{13}}(\Lambda)\nonumber\\
&=&\frac{\partial}{\partial{\Lambda}_{13}}\left\{\int_{-\infty}^{\infty}\int_{-\infty}^t\int_{-\infty}^tf_{\Lambda}(t,t,x_3,x_4)dx_3dx_4dt\right\}\nonumber\\
&=&\int_{-\infty}^{\infty}\int_{-\infty}^t(\partial_{1}f_{\Lambda})(t,t,t,x_4)dx_4dt\nonumber\\
&=&\int_{-\infty}^{\infty}\int_{-\infty}^0(\partial_{1}f_{\Lambda})(t,t,t,x_4+t)dx_4dt\nonumber\\
&=&\int_{-\infty}^0\int_{-\infty}^{\infty}(\partial_{1}f_{\Lambda})(t,t,t,x_4+t)dt dx_4\nonumber\\
&=&\frac{-1}{\sqrt{(2\pi)^{4}{\rm det}({\Lambda})}}\int_{-\infty}^0\int_{-\infty}^{\infty}
\left[\sum_{l=1}^4(\Lambda^{-1})_{1l}t+(\Lambda^{-1})_{14}x_4\right]\nonumber\\
&&\ \ \ \ \ \ \ \cdot e^{-\frac{1}{2}\{\sum_{k,l=1}^4(\Lambda^{-1})_{kl}t^2+2\sum_{l=1}^4(\Lambda^{-1})_{4l}x_4t
+(\Lambda^{-1})_{44}x^2_4\}}dtdx_4\nonumber\\
&=&\frac{-1}{\sqrt{(2\pi)^{4}{\rm det}({\Lambda})}}\int_{-\infty}^0\int_{-\infty}^{\infty}
\left[\sum_{l=1}^4(\Lambda^{-1})_{1l}t+(\Lambda^{-1})_{14}x_4\right]e^{-\frac{1}{2}\left\{\sqrt{\sum_{k,l=1}^4(\Lambda^{-1})_{kl}}t+\frac{
 \sum_{l=1}^4(\Lambda^{-1})_{4l}x_4}{\sqrt{\sum_{k,l=1}^4(\Lambda^{-1})_{kl}}}\right\}^2}dt\nonumber\\
&&\ \ \ \ \ \ \ e^{-\frac{1}{2}\left\{\left[(\Lambda^{-1})_{44}-\frac{\left(\sum_{l=1}^4(\Lambda^{-1})_{4l}\right)^2 }
{{\sum_{k,l=1}^4(\Lambda^{-1})_{kl}}}\right]x^2_4\right\}} dx_4\nonumber\\
&=&\frac{-1}{\sqrt{(2\pi)^{3}{\rm det}({\Lambda})\sum_{k,l=1}^4(\Lambda^{-1})_{kl}}}\int_{-\infty}^0
\left[ \left((\Lambda^{-1})_{14}-\frac{(\sum_{l=1}^4(\Lambda^{-1})_{1l})(\sum_{l=1}^4(\Lambda^{-1})_{4l})}{{\sum_{k,l=1}^4(\Lambda^{-1})_{kl}}}\right)x_4\right]\nonumber\\
&&\ \ \ \ \ \ \cdot
e^{-\frac{1}{2}\left\{\left[(\Lambda^{-1})_{44}-\frac{\left(\sum_{l=1}^4(\Lambda^{-1})_{4l}\right)^2 }
{{\sum_{k,l=1}^4(\Lambda^{-1})_{kl}}}\right]x^2_4\right\}} dx_4\nonumber\\
&=&\frac{1}{\sqrt{(2\pi)^{3}{\rm det}({\Lambda})\sum_{k,l=1}^4(\Lambda^{-1})_{kl}}}
\left[\frac{(\Lambda^{-1})_{14}-\frac{\left(\sum_{l=1}^4(\Lambda^{-1})_{1l}\right)\left(\sum_{l=1}^4(\Lambda^{-1})_{4l}\right)}{{\sum_{k,l=1}^4(\Lambda^{-1})_{kl}}}}
{(\Lambda^{-1})_{44}-\frac{\left(\sum_{l=1}^4(\Lambda^{-1})_{4l}\right)^2 }{{\sum_{k,l=1}^4(\Lambda^{-1})_{kl}}}}\right] \nonumber\\
&=&\frac{(\tilde{\Lambda})_{23}}{2\sqrt{\pi^3}{\tilde A}\{4(\Lambda')_{12}(\Lambda')_{23}
-[(\Lambda')_{12}+(\Lambda')_{23}-(\Lambda')_{13}]^2\}}.
\end{eqnarray*}

By (\ref{ffff}), (\ref{F12}), (\ref{Fou})  and (\ref{dis3}), we get
\begin{eqnarray*}
-\frac{\partial^2 F}{\partial\Lambda_{12}\partial\Lambda_{34}}(\Lambda)&=&\frac{\partial}{\partial{\Lambda}_{34}}\left\{\int_{-\infty}^{\infty}\int_{-\infty}^t\int_{-\infty}^tf_{\Lambda}(t,t,x_3,x_4)dx_3dx_4dt\right\}\nonumber\\
&=&\int_{-\infty}^{\infty}f_{\Lambda}(t,t,t,t)dt\nonumber\\
&=&\frac{1}{\sqrt{(2\pi)^4{\rm det}({\Lambda})}}\int_{-\infty}^{\infty}e^{-\frac{\sum_{k,l=1}^4(\Lambda^{-1})_{kl}}{2}t^2}dt\nonumber\\
&=&\frac{1}{\sqrt{(2\pi)^{3}{\rm det}({\Lambda})\sum_{k,l=1}^4(\Lambda^{-1})_{kl}}}\nonumber\\
&=&\frac{1}{2\sqrt{\pi^3}{\tilde A}}.
\end{eqnarray*}
\hfill\fbox

\vskip 0.2cm
Define
$$
{\cal S}_1=\{\Lambda\in{\bar{\cal S}}:({\Lambda}')_{kl}{\tilde A}^2+(\tilde{\Lambda})^2_{kl}>0\ {\rm for\ all}\ k<l\},
$$
and
$$
{\cal S}_2=\{\Lambda\in{\bar{\cal S}}:{\Lambda}_{kl}\not=1\ {\rm for\ all}\ k<l\}.
$$
\begin{lem}\label{lem0099}
We have
$$
{\cal S}_1={\cal S}_2.
$$
\end{lem}
\noindent {\bf Proof.}\ \ Suppose that $\Lambda\notin {\cal S}_2$. Then, ${\Lambda}_{kl}=1$ for some $k<l$. We assume without loss of generality that ${\Lambda}_{12}=1$.
By (\ref{qqq}), we get $(\tilde{\Lambda})_{12}=0$. Thus, $({\Lambda}')_{12}{\tilde A}^2+(\tilde{\Lambda})^2_{12}=0$, which implies that
$\Lambda\notin {\cal S}_1$.

Suppose that $\Lambda\notin {\cal S}_1$. Then, $({\Lambda}')_{kl}{\tilde A}^2+(\tilde{\Lambda})^2_{kl}=0$ for some $k<l$. We assume without loss of generality that
$({\Lambda}')_{12}{\tilde A}^2+(\tilde{\Lambda})^2_{12}=0$. Take an $X^{\Lambda}\in{\cal G}$ with covariance matrix $\Lambda$. By (\ref{ADD123})--(\ref{deg3}) of the Appendix, we find that
either  $(X^{\Lambda}_1-X^{\Lambda}_2,X^{\Lambda}_3-X^{\Lambda}_2)$ or $(X^{\Lambda}_1-X^{\Lambda}_2,X^{\Lambda}_4-X^{\Lambda}_2)$
is a degenerate 2-dimensional Gaussian random vector. Without loss of generality, we assume
that $(X^{\Lambda}_1-X^{\Lambda}_2,X^{\Lambda}_3-X^{\Lambda}_2)$ is degenerate. Then, there exist constants $\alpha$ and $\beta$ such that $|\alpha|+|\beta|>0$ and
$$
\alpha(X^{\Lambda}_1-X^{\Lambda}_2)+\beta(X^{\Lambda}_3-X^{\Lambda}_2)=0.
$$
If $\beta=0$, then $X^{\Lambda}_1=X^{\Lambda}_2$ and thus ${\Lambda}_{12}=1$; if $\alpha=0$, then $X^{\Lambda}_2=X^{\Lambda}_3$ and
thus ${\Lambda}_{23}=1$. If $\alpha\beta\not=0$ and $\alpha+\beta=0$,
then $X^{\Lambda}_1=X^{\Lambda}_2$ and thus ${\Lambda}_{12}=1$. If $\alpha\beta\not=0$ and $\alpha+\beta\not=0$, we get
\begin{eqnarray}\label{dependent-equation}
X^{\Lambda}_2=\frac{\alpha}{\alpha+\beta}X^{\Lambda}_1+\frac{\beta}{\alpha+\beta}X^{\Lambda}_3.
\end{eqnarray}
By  $E[(X^{\Lambda}_1)^2]=E[(X^{\Lambda}_2)^2]=E[(X^{\Lambda}_3)^2]=1$ and \eqref{dependent-equation}, we get
$$
1=\left(\frac{\alpha}{\alpha+\beta}\right)^2+\left(\frac{\beta}{\alpha+\beta}\right)^2
+\frac{2\alpha\beta}{(\alpha+\beta)^2}{\Lambda}_{13},
$$
which implies that ${\Lambda}_{13}=1$. Hence $\Lambda\notin {\cal S}_2$.\hfill\fbox

\begin{thm}\label{Formula} (i) Suppose that $\Lambda\in{\cal S}_1$. We have
\begin{eqnarray}\label{June1}
F(\Lambda)=\frac{1}{2\sqrt{\pi^3}}\sum_{k<l}\sqrt{({\Lambda}')_{kl}}\arccos
\left(\frac{(\tilde{\Lambda})_{kl}}{\sqrt{({\Lambda}')_{kl}{\tilde A}^2+(\tilde{\Lambda})^2_{kl}}}\right).
\end{eqnarray}
(ii) Suppose that $\Lambda_{14}=1$. We have
\begin{eqnarray}\label{ert0}
F(\Lambda)=\frac{1}{2\sqrt{\pi}}\left[\sqrt{(\Lambda')_{12}}+\sqrt{(\Lambda')_{13}}+\sqrt{(\Lambda')_{23}}\right].
\end{eqnarray}
\end{thm}
\noindent {\bf Proof.}\ \ (i) We first show that (\ref{June1}) holds for $\Lambda\in{\cal S}$. Denote the right hand side of (\ref{June1})
by $\tilde F(\Lambda)$. Let $\Lambda^*$ be the $4\times 4$ matrix with all elements equal to 1. It is easy to see that
 $F(\Lambda^*)={\tilde F}(\Lambda^*)=0$ and $\tilde F$ is continuous at $\Lambda^*$.
 By virtue of characteristic functions, we can show
that $F$ is continuous on ${\bar{\cal S}}$. Since both $F$ and $\tilde F$ are differentiable on ${\cal S}$, to prove that (\ref{June1}) holds on ${\cal S}$,
we need only show that $F$ and $\tilde F$
 have the same partial derivatives on ${\cal S}$. Without loss of generality, in the sequel, we show that
 $$
 \frac{\partial F}{\partial \Lambda_{12}}=\frac{{\partial {\tilde F}}}{\partial \Lambda_{12}}.
 $$
  By (\ref{June0}), we get
$${\tilde F}(\Lambda)=-2\sum_{k<l}({\Lambda}')_{kl}\frac{\partial F}{\partial \Lambda_{kl}}.
$$
Hence we need only show that
\begin{eqnarray}\label{equal2}
\sum_{k<l}({\Lambda}')_{kl}\frac{\partial^2 F}{\partial \Lambda_{12}\partial \Lambda_{kl}}=\frac{1}{2}\frac{\partial F}{\partial \Lambda_{12}},
\ \ \ \ \forall\Lambda\in{\cal S}.
\end{eqnarray}

We use $a,b,c,d,e,f$ to represent $\Lambda_{12}, \Lambda_{13},\Lambda_{14},\Lambda_{23},\Lambda_{24},\Lambda_{34}$, respectively. By (\ref{AA}) and (\ref{ADD123}) of the Appendix, we have that
$$
g:=[4(1-a)(1-d)-(1-a-d+b)^2]\cdot[4(1-a)(1-e)-(1-a-e+c)^2]>0.
$$
By (\ref{June0}) and Theorem \ref{Gauss44}, we find that
 the difference of the two sides of (\ref{equal2}) multiplied by
 $4\sqrt{\pi^3}{\tilde A}g$ is equal to
{\small\begin{eqnarray*}
&&((1-b)^2-(1-b)(1-a+1-c+1-d+1-f-2(1-e))+(a-d)(c-f))(1-a+(1-d)-(1-b))\\
&&\ \ \ \ \cdot(4(1-a)(1-e)-(1-a-e+c)^2)\\
&&+((1-c)^2-(1-c)(1-a+1-b+1-e+1-f-2(1-d))+(a-e)(b-f))(1-a+(1-e)-(1-c))\\
&&\ \ \ \ \cdot (4(1-a)(1-d)-(1-a-d+b)^2)\\
&&-2((1-d)^2-(1-d)(1-a+1-b+1-e+1-f-2(1-c))+(a-b)(e-f))(1-b)\\
&&\ \ \ \ \cdot (4(1-a)(1-e)-(1-a-e+c)^2)\\
&&-2((1-e)^2-(1-e)(1-a+1-c+1-d+1-f-2(1-b))+(a-c)(d-f))(1-c)\\
&&\ \ \ \ \cdot (4(1-a)(1-d)-(1-a-d+b)^2)\\
&&-2((1-b)^2-(1-b)(1-a+1-c+1-d+1-f-2(1-e))+(a-d)(c-f))(1-d)\\
&&\ \ \ \ \cdot (4(1-a)(1-e)-(1-a-e+c)^2)\\
&&-2((1-c)^2-(1-c)(1-a+1-b+1-e+1-f-2(1-d))+(a-e)(b-f))(1-e)\\
&&\ \ \ \ \cdot (4(1-a)(1-d)-(1-a-d+b)^2)\\
&&-2(1-f)(4(1-a)(1-d)-(1-a-d+b)^2)(4(1-a)(1-e)-(1-a-e+c)^2).
\end{eqnarray*} }

\noindent By direct calculation or more quickly by using Mathematica, we can show that the above multivariate polynomial of $(a,b,\dots,f)$ is equal to 0. We refer the reader to \S 7.3 for the Mathematica code. Then,  (\ref{equal2}) holds and thus (\ref{June1}) holds for $\Lambda\in{\cal S}$.

Since $F$ is continuous on ${\bar{\cal S}}$ and $\tilde F$ is continuous on ${\cal S}_1$. By the density of ${\cal S}$ in $\bar{\cal S}$, we conclude that
(\ref{June1}) holds for all $\Lambda\in{\cal S}_1$.

(ii) Suppose that ${\Lambda}_{14}=1$. We have
\begin{eqnarray}\label{max11}
F(\Lambda)&=&E[\max(X^{\Lambda}_1,X^{\Lambda}_2,X^{\Lambda}_3,X^{\Lambda}_4)]\nonumber\\
&=&E[\max(X^{\Lambda}_1,X^{\Lambda}_2,X^{\Lambda}_3)]\nonumber\\
&=&E[\max(X^{\Lambda}_1-X^{\Lambda}_3,X^{\Lambda}_2-X^{\Lambda}_3,0)]\nonumber\\
&=&\frac{E[\max(X^{\Lambda}_1-X^{\Lambda}_3,X^{\Lambda}_2-X^{\Lambda}_3)]+E[|\max(X^{\Lambda}_1-X^{\Lambda}_3,X^{\Lambda}_2-X^{\Lambda}_3)|]}{2},
\end{eqnarray}
and
\begin{eqnarray}\label{max12}
&&E[\max(X^{\Lambda}_1-X^{\Lambda}_3,X^{\Lambda}_2-X^{\Lambda}_3)]\nonumber\\
&=&\frac{E[(X^{\Lambda}_1-X^{\Lambda}_3)+(X^{\Lambda}_2-X^{\Lambda}_3)]+E[|(X^{\Lambda}_1-X^{\Lambda}_3)-(X^{\Lambda}_2-X^{\Lambda}_3)|]}{2}\nonumber\\
&=&\frac{E[|X^{\Lambda}_1-X^{\Lambda}_2|]}{2}\nonumber\\
&=&\sqrt{\frac{(\Lambda')_{12}}{\pi}}.
\end{eqnarray}
Note that $(X^{\Lambda}_1-X^{\Lambda}_3,X^{\Lambda}_2-X^{\Lambda}_3)$ and $(X^{\Lambda}_3-X^{\Lambda}_1,X^{\Lambda}_3-X^{\Lambda}_2)$
have the same distribution. Then,
\begin{eqnarray}\label{max13}
&&E[|\max(X^{\Lambda}_1-X^{\Lambda}_3,X^{\Lambda}_2-X^{\Lambda}_3)|]\nonumber\\
&=&\frac{E[|\max(X^{\Lambda}_1-X^{\Lambda}_3,X^{\Lambda}_2-X^{\Lambda}_3)|]+E[|\max(X^{\Lambda}_3-X^{\Lambda}_1,X^{\Lambda}_3-X^{\Lambda}_2)|]}{2}\nonumber\\
&=&\frac{E[|\max(X^{\Lambda}_1-X^{\Lambda}_3,X^{\Lambda}_2-X^{\Lambda}_3)|]+E[|\min(X^{\Lambda}_1-X^{\Lambda}_3,X^{\Lambda}_2-X^{\Lambda}_3)|]}{2}\nonumber\\
&=&\frac{E[|X^{\Lambda}_1-X^{\Lambda}_3|]+E[|X^{\Lambda}_2-X^{\Lambda}_3|]}{2}\nonumber\\
&=&\sqrt{\frac{(\Lambda')_{13}}{\pi}}+\sqrt{\frac{(\Lambda')_{23}}{\pi}}.
\end{eqnarray}
Therefore, (\ref{ert0}) holds by (\ref{max11})--(\ref{max13}).\hfill\fbox

\begin{rem}\label{rempo}
Note that if $\Lambda\in\bar{\cal S}$ satisfying $\Lambda_{kl}=r$ for all $k<l$ and some $r\in \mathbb{R}$, then $-\frac{1}{3}\le r\le 1$. If $-\frac{1}{3}\le r< 1$, by (\ref{June1}), (\ref{mmmm})--(\ref{ADD123}) of the Appendix, and elementary calculation, we get
\begin{equation}\label{add09}
F(\Lambda)=\frac{3\sqrt{1-r}\arccos\left(-\frac{1}{3}\right)}{\sqrt{\pi^3}}.
\end{equation}
Obviously, (\ref{add09}) also holds for the case that $r=1$.
\end{rem}

\vskip 0.3cm
Denote by ${\cal G}_3$ the set of all three-dimensional centered
Gaussian random vectors $X=(X_1,X_2,X_3)$ with $E[X_1^2]=E[X_2^2]=E[X_3^2]=1$. Define
$$
{\cal M}(X)=E[\max(X_1,X_2,X_3)], \ \ \ \ X\in {\cal G}_3.
$$
We use $\bar{\cal S}_3$ to denote the set of all nonnegative-definite
 $3\times 3$ symmetric real
matrices with all diagonal elements equal to 1. For $\Lambda\in{\bar{\cal S}}_3$,
take an $X^{\Lambda}\in{\cal G}_3$ with covariance matrix $\Lambda$ and define
$$
F(\Lambda)={\cal M}(X^{\Lambda}).
$$
\begin{cor}\label{cor3}
The function ${\cal M}$ defined on ${\cal G}_3$ attains its maximum at $X$ if and only if
$\Lambda^X_{kl}=-\frac{1}{2}$ for all $k<l$.
\end{cor}

\noindent {\bf Proof.}\ \ By (\ref{ert0}), we know that $F(\Lambda)$ is a strictly concave function defined on $\bar{\cal S}_3$. Let $\Lambda\in \bar{\cal S}_3$. Define ${\bar\Lambda}\in \bar{\cal S}_3$ by
\begin{eqnarray*}
{\bar\Lambda}_{kl}=\left\{
  \begin{array}{ll}
1, \ \ \ \ & {\rm if}\ \ k=l, \\
\frac{\sum_{k<l}\Lambda_{kl}}{3}, \ \ \ \ \ \ & {\rm otherwise}.
  \end{array}
\right.
\end{eqnarray*}
Then, we have that $F(\Lambda)\le F({\bar\Lambda})$ and the equality holds if and only if $\Lambda={\bar\Lambda}$.
Note that if $\Lambda\in\bar{\cal S}_3$ satisfying $\Lambda_{kl}=r$ for all $k<l$ and some $r\in \mathbb{R}$, then $-\frac{1}{2}\le r\le 1$. Therefore, the proof is complete by (\ref{ert0}).\hfill\fbox

\section{Proof of Theorem \ref{main}}\setcounter{equation}{0}

\noindent {\it Step 1.}\ \ First, we show that if ${\cal M}(X)$ attains its maximum at $X\in{\cal G}$ then $X$ is degenerate, i.e., $X_1,X_2,X_3,X_4$ are linear dependent.
In fact, if $X$ is non-degenerate then $\Lambda^X\in{\cal S}$. To simplify notation, in the sequel, we denote $\Lambda^X$ by $\Lambda$. For $\varepsilon>0$, define $\Lambda^{\varepsilon}$ by
\begin{eqnarray*}
\Lambda^{\varepsilon}_{kl}=\left\{
  \begin{array}{ll}
\Lambda_{12}-\varepsilon, \ \ \ \ & {\rm if}\ \ (k,l)=(1,2)\ {\rm or}\ (2,1), \\
\Lambda_{kl}, \ \ \ \ \ \ & {\rm otherwise}.
  \end{array}
\right.
\end{eqnarray*}
If $\varepsilon$ is small enough, then $\Lambda^{\varepsilon}\in{\cal S}$.  By Lemma \ref{Gauss2},
we get $F(\Lambda)<F(\Lambda^{\varepsilon})$ and hence  arrive at a contradiction.

Next, we show that if ${\cal M}(X)$ attains its maximum at $X\in{\cal G}$, then $\Lambda_{kl}\not=1$ for any $k<l$. Assume that, for example, $\Lambda_{14}=1$. By Corollary \ref{cor3}, we get
$$
F(\Lambda)\le \frac{3\sqrt{\frac{3}{2}}}{2\sqrt{\pi}}
\approx  1.036482< 1.18862\approx\frac{3\sqrt{1+\frac{1}{3}}\arccos(-\frac{1}{3})}{\sqrt{\pi^3}}.
$$
Hence $X$ is not a maximum point by Remark \ref{rempo} and we arrive at a contradiction. Thus, $\Lambda\in {\cal S}_1$ by Lemma \ref{lem0099}.

We regard
$X_1,X_2,X_3,X_4$ as vectors of the space $L^2(\Omega,{\cal F},P)$. Denote  by ${\cal F}_1,{\cal F}_2,{\cal F}_3,{\cal F}_4$ the facets of the tetrahedron $\{X_1,X_2,X_3,X_4\}$ with vertexes $\{X_1,X_2,X_3\}$, $\{X_1,X_3,X_4\}$, $\{X_1,X_2,X_4\}$, $\{X_2,X_3,X_4\}$, respectively. We use $\alpha_{ij}$, $i<j$, to denote  the {outer} dihedral angle of ${\cal F}_i$ and ${\cal F}_j$.
 By (\ref{qqq}) and (\ref{ADD123}) of the Appendix, we get
{\small\begin{eqnarray}\label{geometric-expression}
&&\frac{(\tilde{\Lambda})_{12}}{\sqrt{({\Lambda}')_{12}{\tilde A}^2+(\tilde{\Lambda})^2_{12}}}\nonumber\\
&=&\frac{(1+\Lambda_{13}-\Lambda_{12}-\Lambda_{23})(1+\Lambda_{14}-\Lambda_{12}-\Lambda_{24})-2(1-\Lambda_{12})
(1+\Lambda_{34}-\Lambda_{23}-\Lambda_{24})}{\sqrt{
\left[4(1-\Lambda_{12})(1-\Lambda_{23})-\left(1+\Lambda_{13}-\Lambda_{12}-\Lambda_{23}\right)^2\right]\left[4(1-\Lambda_{12})(1-\Lambda_{24})-\left(1+\Lambda_{14}-\Lambda_{12}-\Lambda_{24}\right)^2\right]
}}\nonumber\\
&=&\{{\rm Cov}(X_1-X_2,X_3-X_2){\rm Cov}(X_1-X_2,X_4-X_2)-{\rm Var}(X_1-X_2){\rm Cov}(X_3-X_2,X_4-X_2)\}\nonumber\\
&&\cdot\{{\rm Var}(X_1-X_2){\rm Var}(X_3-X_2)-[{\rm Cov}(X_1-X_2,X_3-X_2)]^2\}^{-\frac{1}{2}}\nonumber\\
&&\cdot\{{\rm Var}(X_1-X_2){\rm Var}(X_4-X_2)-[{\rm Cov}(X_1-X_2,X_4-X_2)]^2\}^{-\frac{1}{2}}.
\end{eqnarray}}
\vskip -0.3cm
\noindent Denote by $\beta_1,\beta_2,\beta_3$ the angles between $\overrightarrow{X_2X_1}$ and $\overrightarrow{X_2X_3}$, $\overrightarrow{X_2X_1}$ and $\overrightarrow{X_2X_4}$, $\overrightarrow{X_2X_3}$ and $\overrightarrow{X_2X_4}$, respectively. Then, we obtain by
 \eqref{geometric-expression}  that
\begin{equation}\label{GEOMP}
\frac{(\tilde{\Lambda})_{12}}{\sqrt{({\Lambda}')_{12}{\tilde A}^2+(\tilde{\Lambda})^2_{12}}}
=\frac{\cos\beta_1\cos\beta_2-\cos \beta_3}{\sin\beta_1\sin\beta_2}=\cos\alpha_{13}.
\end{equation}
Similarly, we can show that
\begin{eqnarray}\label{GEOM1}
\frac{(\tilde{\Lambda})_{13}}{\sqrt{({\Lambda}')_{13}{\tilde A}^2+(\tilde{\Lambda})^2_{13}}}
&=&\cos\alpha_{12},\nonumber\\
\frac{(\tilde{\Lambda})_{14}}{\sqrt{({\Lambda}')_{14}{\tilde A}^2+(\tilde{\Lambda})^2_{14}}}
&=&\cos\alpha_{23},\nonumber\\
\frac{(\tilde{\Lambda})_{23}}{\sqrt{({\Lambda}')_{23}{\tilde A}^2+(\tilde{\Lambda})^2_{23}}}
&=&\cos\alpha_{14},\nonumber\\
\frac{(\tilde{\Lambda})_{24}}{\sqrt{({\Lambda}')_{24}{\tilde A}^2+(\tilde{\Lambda})^2_{24}}}
&=&\cos\alpha_{34},\nonumber\\
\frac{(\tilde{\Lambda})_{34}}{\sqrt{({\Lambda}')_{34}{\tilde A}^2+(\tilde{\Lambda})^2_{34}}}
&=&\cos\alpha_{24}.
\end{eqnarray}

Now we show that if ${\cal M}(X)$ attains its maximum at $X\in{\cal G}$, then $X_1,X_2,X_3,X_4$ do not lie in the same plane. Assume that $X_1,X_2,X_3,X_4$ lie in the same plane. Then, two outer  dihedral angles of the degenerate tetrahedron $\{X_1,X_2,X_3,X_4\}$ equal 0. Hence, we obtain by
Theorem \ref{Formula} that
$$
F(\Lambda^X)\le\frac{4\pi}{2\sqrt{\pi^3}}
\approx 1.128379< 1.18862\approx\frac{3\sqrt{1+\frac{1}{3}}\arccos(-\frac{1}{3})}{\sqrt{\pi^3}},
$$
which implies that $X$ is not a maximum point  by Remark \ref{rempo} and we arrive at a contradiction.

\noindent {\it Step 2.}\ \  Suppose that ${\cal M}(X)$ attains its maximum at $X\in{\cal G}$. Based on the analysis of Step 1, we know that $X_1,X_2,X_3,X_4$ are linear dependent and do not lie in the same plane.

First, we focus on the vertex $X_1$. We consider two different cases. Case 1: $X_i\not=-X_1$ for any $2\le i\le 4$.
Denote by ${\cal P}_1$, ${\cal P}_2$ and ${\cal P}_3$ the planes determined by $\{0,X_1,X_2\}$, $\{0,X_1,X_3\}$ and $\{0,X_1,X_4\}$, respectively. At least one of the three planes ${\cal P}_1$, ${\cal P}_2$, ${\cal P}_3$ is different from the other two planes. We assume without loss of generality that ${\cal P}_1\not={\cal P}_2$ and ${\cal P}_1\not={\cal P}_3$. Case 2: $X_i=-X_1$ for some $2\le i\le 4$, say $X_4=-X_1$. Denote by ${\cal P}_1$ the plane determined by $\{0,X_1,X_2\}$ and  by  ${\cal P}_2={\cal P}_3$ the plane determined by $\{0,X_1,X_3,X_4\}$, respectively.

Denote $\theta=(\theta_1,\theta_2,\theta_3,\theta_4,\theta_5)$. Take three independent standard Gaussian random variables $U_1,U_2,U_3$ as the basis and  set  $(X_1,X_2,X_3,X_4)$ to be
\begin{eqnarray*}
X_1(\theta)&=&(0,0,1)^T,\nonumber\\
X_2(\theta)&=&(\sin\theta_1,0,\cos\theta_1)^T,\nonumber\\
X_3(\theta)&=&(\sin\theta_2\cos\theta_4,\sin\theta_2\sin\theta_4,\cos\theta_2)^T,\nonumber\\
X_4(\theta)&=&(\sin\theta_3\cos\theta_5,\sin\theta_3\sin\theta_5,\cos\theta_3)^T,
\end{eqnarray*}
where $\theta_1\in(0,\pi)$ denotes the angle between $\overrightarrow{0X_1}$ and $\overrightarrow{0X_2}$, $\theta_2,\theta_3\in(0,2\pi)$ denote the anti-clockwise angles from $\overrightarrow{0X_1}$ to $\overrightarrow{0X_3}$ and from $\overrightarrow{0X_1}$ to $\overrightarrow{0X_4}$, respectively, and $\theta_4,\theta_5\in(0,\pi)$ denote the anti-clockwise angles from ${\cal P}_1$ to ${\cal P}_2$ and from ${\cal P}_1$ to ${\cal P}_3$, respectively.
Then,
\begin{eqnarray*}
\Lambda_{12}(\theta)&=&\cos\theta_1,\nonumber\\
\Lambda_{13}(\theta)&=&\cos\theta_2,\nonumber\\
\Lambda_{14}(\theta)&=&\cos\theta_3,\nonumber\\
\Lambda_{23}(\theta)&=&\sin\theta_1\sin\theta_2\cos\theta_4+\cos\theta_1\cos\theta_2,\nonumber\\
\Lambda_{24}(\theta)&=&\sin\theta_1\sin\theta_3\cos\theta_5+\cos\theta_1\cos\theta_3,\nonumber\\
\Lambda_{34}(\theta)&=&\sin\theta_2\sin\theta_3\cos(\theta_4-\theta_5)+\cos\theta_2\cos\theta_3.
\end{eqnarray*}

For $T\in{\cal S}_1$ and $k<l$, define
\begin{equation}\label{gai2}
F_{kl}(T)
=-\frac{1}{4\sqrt{\pi^3(T')_{kl}}}\arccos\left(\frac{(\tilde{T})_{kl}}{\sqrt{({T}')_{kl}{\tilde A}^2+(\tilde{T})^2_{kl}}}\right).
\end{equation}
Then, $F_{kl}$ is a continuous function on ${\cal S}_1$. By (\ref{June0}), we get
$$
F_{kl}(T)=\frac{\partial F}{\partial \Lambda_{kl}}(T),\ \ \ \ \forall T\in{\cal S}.
$$

For $T\in{\bar{\cal S}}$ and $0<\delta<1$, define $T^{\delta}=(1-\delta)T+\delta I$, where $I$ is the identity matrix. Then, $T^{\delta}\in{\cal S}$ and $T^{\delta}$ approximates $T$ as $\delta\downarrow0$. For $1\le i\le 5$, denote by $e_i$ the five-dimensional unit vector with the $i$-th component equal to 1 and the other four components equal to 0. Define
$$
R(\theta)=F(\Lambda(\theta)).
$$
By the continuity of the functions $F_{kl}$, $k<l$, we get
\begin{eqnarray*}
&&\frac{\partial R}{\partial \theta_i}\nonumber\\
&=&\lim_{\varepsilon\rightarrow 0}\frac{F[\Lambda(\theta+\varepsilon e_i)]-F[\Lambda(\theta)]}{\varepsilon}\nonumber\\
&=&\lim_{\varepsilon\rightarrow 0}\lim_{\delta\rightarrow 0}\frac{F[\Lambda^{\delta}(\theta+\varepsilon e_i)]-F[\Lambda^{\delta}(\theta)]}{\varepsilon}\nonumber\\
&=&\lim_{\varepsilon\rightarrow 0}\lim_{\delta\rightarrow 0}\frac{1}{\varepsilon}
\int_0^1d{F[p\Lambda^{\delta}(\theta+\varepsilon e_i)+(1-p)\Lambda^{\delta}(\theta)]}\nonumber\\
&=&\sum_{k<l}\left\{\lim_{\varepsilon\rightarrow 0}\lim_{\delta\rightarrow 0}
\left(\int_0^1\frac{\partial F}{\partial\Lambda_{kl}}{[p\Lambda^{\delta}(\theta+\varepsilon e_i)+(1-p)\Lambda^{\delta}(\theta)]\cdot \frac{(1-\delta)[\Lambda_{kl}(\theta+\varepsilon e_i)-\Lambda_{kl}(\theta)]}{\varepsilon}}dp\right)\right\}\nonumber\\
&=&\sum_{k<l}\left\{\lim_{\varepsilon\rightarrow 0}\lim_{\delta\rightarrow 0}
\int_0^1\frac{\partial F}{\partial\Lambda_{kl}}{[p\Lambda^{\delta}(\theta+\varepsilon e_i)+(1-p)\Lambda^{\delta}(\theta)]dp\cdot\lim_{\varepsilon\rightarrow 0}\lim_{\delta\rightarrow 0} \frac{(1-\delta)[\Lambda_{kl}(\theta+\varepsilon e_i)-\Lambda_{kl}(\theta)]}{\varepsilon}}\right\}\nonumber\\
&=&\sum_{k<l}F_{kl}\frac{\partial \Lambda_{kl}}{\partial \theta_i}.
\end{eqnarray*}
Since ${\cal M}(X)$ attains its maximum at $X$, we get
$$
\frac{\partial R}{\partial \theta_i}=0, \ \ \ \ i=4,5,
$$
which implies that
{\begin{eqnarray}\label{derivative-4}
&&F_{23}\sin\theta_1\sin\theta_2\sin\theta_4
+F_{34}\sin\theta_2\sin\theta_3\sin(\theta_4-\theta_5)=0,\nonumber\\
&&F_{24}\sin\theta_1\sin\theta_3\sin\theta_5-F_{34}\sin\theta_2\sin\theta_3\sin(\theta_4-\theta_5)=0.
\end{eqnarray}}

For $1\le i\le 4$, let $0P_i$ be the perpendicular line from 0 to ${\cal F}_i$ with
$P_i\in {\cal F}_i$ and denote its length by
${\cal L}_i$. Note that, for $i<j$, the angle between $0P_i$ and $0P_j$ equals $\alpha_{ij}$, the {outer} dihedral angle of ${\cal F}_i$ and ${\cal F}_j$.
Denote by $V_1$, $V_2$, $V_3$ and $V_4$ the volumes of the
tetrahedrons with vertexes $\{0,X_1,X_2,X_3\}$, $\{0,X_1,X_3,X_4\}$, $\{0,X_1,X_2,X_4\}$ and $\{0,X_2,X_3,X_4\}$, respectively. We have
\begin{eqnarray}\label{haha3}
V_1&=&\frac{1}{6}\sin\theta_1\cdot|\sin\theta_2|\cdot\sin\theta_4,\nonumber\\
V_2&=&\frac{1}{6}|\sin\theta_2\sin\theta_3\sin(\theta_4-\theta_5)|,\nonumber\\
V_3&=&\frac{1}{6}\sin\theta_1\cdot|\sin\theta_3|\cdot\sin\theta_5.
\end{eqnarray}
Since each $F_{kl}$ is strictly negative, we obtain by (\ref{derivative-4}) and (\ref{haha3}) that
{\begin{eqnarray}\label{supp1}
F_{23}V_1=
F_{34}V_2=F_{24}V_3.
\end{eqnarray}}
\vskip -0.5cm
\noindent This implies that $V_1V_2V_3\not=0$ since $V_1>0$.

Next, by focusing on the vertexes $X_2,X_3,X_4$, respectively, and following the above argument, we obtain that
{\begin{eqnarray}\label{com2}
&&F_{14}V_3= F_{34}V_4=F_{13}V_1,\nonumber\\
&&F_{24}V_4=F_{12}V_1=
F_{14}V_2,\nonumber\\
&&F_{13}V_2=F_{12}V_3=
F_{23}V_4,
\end{eqnarray}}
and $V_1V_2V_3V_4\not=0$. Combining (\ref{supp1}) and (\ref{com2}), we get
\begin{eqnarray}\label{supp2}
&&V_2=V_1\cdot\frac{F_{12}}{F_{14}}=V_1\cdot\frac{F_{23}}{F_{34}},\nonumber\\
&&V_3=V_1\cdot\frac{F_{13}}{F_{14}}=V_1\cdot\frac{F_{23}}{F_{24}},\nonumber\\
&&V_4=V_1\cdot\frac{F_{12}}{F_{24}}
=V_1\cdot\frac{F_{13}}{F_{34}}.
\end{eqnarray}
Further, we obtain by  (\ref{supp2}) that
\begin{eqnarray}
F_{13}&=&F_{12}\cdot\frac{V_3}{V_2},\label{tele1}\\
F_{14}&=&F_{12}\cdot\frac{V_1}{V_2},\label{tele2}\\
F_{23}&=&F_{12}\cdot\frac{V_3}{V_4},\label{tele3}\\
F_{24}&=&F_{12}\cdot\frac{V_1}{V_4},\label{tele4}\\
F_{34}&=&F_{12}\cdot\frac{V_1V_3}{V_2
V_4}.\label{tele5}
\end{eqnarray}

By (\ref{tele1}) and (\ref{tele2}), we get
\begin{equation}\label{K34}
\frac{F_{12}}{V_2}=\frac{F_{13}}{V_3}=
\frac{F_{14}}{V_1}.
\end{equation}
By the law of sines for tetrahedrons, we have that
\begin{equation}\label{K35}
\frac{|\overrightarrow{X_1X_2}|\cdot{\rm Area}({\cal F}_2)}{\sin\alpha_{13}}=\frac{|\overrightarrow{X_1X_3}|\cdot{\rm Area}({\cal F}_3)}{\sin\alpha_{12}}=\frac{|\overrightarrow{X_1X_4}|\cdot{\rm Area}({\cal F}_1)}{\sin\alpha_{23}},
\end{equation}
where $|\cdot|$ denotes the length of a vector. Then, we obtain by (\ref{GEOMP})--(\ref{gai2}), (\ref{K34}) and (\ref{K35}) that
\begin{eqnarray}\label{symmm}
\frac{\alpha_{13}}{{\cal L}_2\sin\alpha_{13}}=\frac{\alpha_{12}}{{\cal L}_3\sin\alpha_{12}}=\frac{\alpha_{23}}{{\cal L}_1\sin\alpha_{23}}.
\end{eqnarray}
Further, by (\ref{tele1})--(\ref{tele5}) and symmetry, we get
\begin{eqnarray*}
&&\frac{\alpha_{13}}{{\cal L}_4\sin\alpha_{13}}=\frac{\alpha_{14}}{{\cal L}_3\sin\alpha_{14}}=\frac{\alpha_{34}}{{\cal L}_1\sin\alpha_{34}},\\
&&\frac{\alpha_{12}}{{\cal L}_4\sin\alpha_{12}}=\frac{\alpha_{14}}{{\cal L}_2\sin\alpha_{14}}=\frac{\alpha_{24}}{{\cal L}_1\sin\alpha_{24}},\\
&&\frac{\alpha_{23}}{{\cal L}_4\sin\alpha_{23}}=\frac{\alpha_{34}}{{\cal L}_2\sin\alpha_{34}}=\frac{\alpha_{24}}{{\cal L}_3\sin\alpha_{24}},
\end{eqnarray*}
which together with (\ref{symmm}) implies that
\begin{eqnarray}\label{WEI}
\frac{{\cal L}_1{\cal L}_2\alpha_{12}}{\sin\alpha_{12}}=\frac{{\cal L}_1{\cal L}_3\alpha_{13}}{\sin\alpha_{13}}=\frac{{\cal L}_1{\cal L}_4\alpha_{14}}{\sin\alpha_{14}}=
\frac{{\cal L}_2{\cal L}_3\alpha_{23}}{\sin\alpha_{23}}=\frac{{\cal L}_2{\cal L}_4\alpha_{24}}{\sin\alpha_{24}}=\frac{{\cal L}_3{\cal L}_4\alpha_{34}}{\sin\alpha_{34}}.
\end{eqnarray}

\noindent {\it Step 3.}\ \ Note that $X_1$ satisfies the following equations:
$$
\langle X_1-P_i,P_i\rangle  =0,\ \ \ \ 1\le i\le 3.
$$
Then, we have that
\begin{equation}\label{LAN1}
\langle X_1,P_i\rangle  =\|P_i\|^2={\cal L}_i^2,\ \ \ \ 1\le i\le 3.
\end{equation}

Define the $3\times 3$ matrix $\Gamma$ by
$$
\Gamma=(P_1,P_2,P_3).
$$
Given a matrix $M$, we denote by ${\rm diag}(M)$ the row vector consisting of the diagonal elements of $M$. Then,  (\ref{LAN1}) becomes
$$
X^T_1\Gamma={\rm diag}(\Gamma^T\Gamma),
$$
which implies that
$$
X^T_1={\rm diag}(\Gamma^T\Gamma)\cdot \Gamma^{-1}.
$$
Define
$$
\Sigma=\Gamma^T\Gamma.
$$
Then, we obtain by $\|X_1\|=1$ that
\begin{equation}\label{LAN2}
{\rm diag}(\Sigma)\cdot\Sigma^{-1}\cdot [{\rm diag}(\Sigma)]^T=1.
\end{equation}

Define
$$
\tilde{P}_i=\frac{P_i}{{\cal L}_i},\ \ 1\le i\le 3,\ \ \ \ {\tilde\Gamma}=({\tilde P}_1,{\tilde P}_2,{\tilde P}_3),\ \ \ \ {\tilde\Sigma}={\tilde\Gamma}^T{\tilde\Gamma},\ \ \ \ \Delta=\left(\begin{array}{ccc}
{\cal L}_1&0&0\\
0&{\cal L}_2&0\\
0&0&{\cal L}_3
\end{array}\right).
$$
We have
$$
\Sigma=\Delta{\tilde\Sigma}\Delta,\ \ \ \ \Sigma^{-1}=\Delta^{-1}{\tilde\Sigma}^{-1}\Delta^{-1}.
$$
Note that
$$
{\rm diag}(\Delta M)={\rm diag}( M\Delta)={\rm diag}(M)\Delta.
$$
Then,
$$
{\rm diag}(\Sigma)={\rm diag}({\tilde\Sigma})\Delta^2=(1,1,1)\Delta^2.
$$
Hence, we obtain by (\ref{LAN2}) that
\begin{eqnarray}\label{LAN3}
(1,1,1)\Delta({\tilde\Sigma})^{-1}\Delta\left(\begin{array}{c}
1\\
1\\
1\end{array}\right)=1.
\end{eqnarray}
Note that
\begin{eqnarray*}
{\tilde\Sigma}=\left(\begin{array}{ccc}
1&\langle {\tilde P}_1, {\tilde P}_2\rangle&\langle {\tilde P}_1, {\tilde P}_3\rangle\\
\langle {\tilde P}_1, {\tilde P}_2\rangle&1&\langle {\tilde P}_2, {\tilde P}_2\rangle\\
\langle {\tilde P}_1, {\tilde P}_3\rangle&\langle {\tilde P}_2, {\tilde P}_2\rangle&1
\end{array}\right)=\left(\begin{array}{ccc}
1&\cos\alpha_{12}&\cos\alpha_{13}\\
\cos\alpha_{12}&1&\cos\alpha_{23}\\
\cos\alpha_{13}&\cos\alpha_{23}&1
\end{array}\right),
\end{eqnarray*}
which is a positive-definite matrix. Thus, we obtain by (\ref{LAN3}) that
\begin{eqnarray}\label{RR11}
({\cal L}_1,{\cal L}_2,{\cal L}_3)\left(\begin{array}{ccc}
1&\cos\alpha_{12}&\cos\alpha_{13}\\
\cos\alpha_{12}&1&\cos\alpha_{23}\\
\cos\alpha_{13}&\cos\alpha_{23}&1
\end{array}\right)^{-1}\left(\begin{array}{c}
{\cal L}_1\\
{\cal L}_2\\
{\cal L}_3\end{array}\right)=1.
\end{eqnarray}

From now on till the end of this paper, we adopt the convention that $\frac{0}{0}=1$ and an increasing (decreasing) function means a strictly increasing (decreasing) function. Define
\begin{equation}\label{422}
f(x)=\frac{\sqrt{1-x^2}}{\arccos x}, \ \ \ \ -1\le x\le1.
\end{equation}
Then, $f$ is an increasing function on $[-1,1]$ with $f(-1)=0$ and $f(1)=1$. Let $f^{-1}$ be the inverse function of $f$. Then, $f^{-1}$
is an increasing function on $[0,1]$ with $f^{-1}(0)=-1$ and $f^{-1}(1)=1$.

For $x,y,z\in (0,\infty)$ satisfying $xy<1,xz<1,yz<1$, define
\begin{eqnarray}\label{423}
\Gamma(x,y,z)=\left(\begin{array}{ccc}
1&f^{-1}(xy)&f^{-1}(xz)\\
f^{-1}(xy)&1&f^{-1}(yz)\\
f^{-1}(xz)&f^{-1}(yz)&1
\end{array}\right).
\end{eqnarray}
Further, if ${\rm det}(\Gamma(x,y,z))>0$, define
\begin{eqnarray}\label{WW5}
H(x,y,z)=(x,y,z)[\Gamma(x,y,z)]^{-1}\left(\begin{array}{c}
x\\
y\\
z
\end{array}\right).
\end{eqnarray}

Denote
$$
r_{ij}=\cos\alpha_{ij},\ \ \ \ i<j.
$$
Then,
$$
f(r_{ij})=\frac{\sin\alpha_{ij}}{\alpha_{ij}}.
$$
By (\ref{WEI}), we may set
\begin{equation}\label{WEI11}
\gamma:=\frac{{\cal L}_1{\cal L}_2\alpha_{12}}{\sin\alpha_{12}}=\frac{{\cal L}_1{\cal L}_3\alpha_{13}}{\sin\alpha_{13}}=\frac{{\cal L}_1{\cal L}_4\alpha_{14}}{\sin\alpha_{14}}=
\frac{{\cal L}_2{\cal L}_3\alpha_{23}}{\sin\alpha_{23}}=\frac{{\cal L}_2{\cal L}_4\alpha_{24}}{\sin\alpha_{24}}=\frac{{\cal L}_3{\cal L}_4\alpha_{34}}{\sin\alpha_{34}}.
\end{equation}
Then, we have that
\begin{equation}\label{WEI2}
{\cal L}_i{\cal L}_j=\gamma f(r_{ij}),\ \ \ \ i<j.
\end{equation}

Define
\begin{equation}\label{WEI1111}
u_i=\frac{{\cal L}_i}{\sqrt{\gamma}},\ \ \ \ 1\le i\le 4.
\end{equation}
By (\ref{WEI2}), we get
\begin{equation}\label{WEI3}
u_iu_j=f(r_{ij}),\ \ r_{ij}=f^{-1}(u_iu_j),\ \ \ \ i<j.
\end{equation}
Then, we obtain by (\ref{RR11}) and (\ref{WEI3}) that
$$
H(u_1,u_2,u_3)=\frac{1}{\gamma}.
$$
Further, by symmetry, we can show that
\begin{eqnarray}\label{HHHH}
H(u_1,u_2,u_3)=H(u_1,u_2,u_4)=H(u_1,u_3,u_4)=H(u_2,u_3,u_4)=\frac{1}{\gamma}.
\end{eqnarray}

\noindent {\it Step 4.}\ \ By (\ref{HHHH}) and Theorem  \ref{maint} of the next section, we conclude that $$u_1=u_2=u_3=u_4,$$
which together with (\ref{WEI11})--(\ref{WEI1111}) implies that
\begin{equation}\label{VVV}
\alpha_{12}=\alpha_{13}=\cdots=\alpha_{34}.
\end{equation}
It is well-known that if the six dihedral angles of a given tetrahedron are congruent, then the tetrahedron is regular (cf. e.g. \cite{USAMO}). Then,
we obtain by (\ref{VVV}) that
$$
(\Lambda')_{12}=(\Lambda')_{13}=\cdots=(\Lambda')_{34}.
$$
Therefore,  $\Lambda_{kl}=-\frac{1}{3}$ for all $k<l$ by Remark \ref{rempo}.\hfill\fbox

\section{Monotonicity of $H$}\setcounter{equation}{0}

Let the function $H$ be defined as in (\ref{WW5}).
For $x,y\in (0,\infty)$ satisfying $xy<1$, define
$$
U_{xy}=\left\{z\in\left(0,\frac{1}{\max\{x,y\}}\right): {\rm det}(\Gamma(x,y,z))>0\right\}.
$$
In this section, we will establish the following crucial result, which has been used in the above section to prove Theorem \ref{main}.

\begin{thm}\label{maint} Let $w_1,w_2\in (0,\infty)$ satisfying $w_1w_2<1$. Suppose that $w_3,w_4\in U_{w_1w_2}$ satisfying $w_3<w_4$. Then, $H(w_1,w_2,w_3)>H(w_1,w_2,w_4)$.
\end{thm}

\noindent {\bf Proof.}\ \ Let $(x,y,z)\in(0,\infty)^3$ satisfying $xy,xz,yz<1$. Define
\begin{equation}\label{51}
\varsigma=f^{-1}(xy), \ \ \ \ \eta=f^{-1}(xz),\ \ \ \ \xi=f^{-1}(yz),
\end{equation}
and
\begin{equation}\label{52}
\theta=\arccos \varsigma, \ \ \ \ \nu=\arccos\eta,\ \ \ \ \mu=\arccos\xi.
\end{equation}
Then, $(\theta,\nu,\mu)\in(0,\pi)^3$ and
\begin{equation}\label{determine}
x=\sqrt{\frac{f(\cos\theta)f(\cos\nu)}{f(\cos\mu)}},\ \ \ \ y=\sqrt{\frac{f(\cos\theta)f(\cos\mu)}{f(\cos\nu)}},\ \ \ \ z=\sqrt{\frac{f(\cos\mu)f(\cos\nu)}{f(\cos\theta)}}.
\end{equation}
For $(\theta,\nu,\mu)\in(0,\pi)^3$, we define
$$
J(\theta,\nu,\mu)=H(x,y,z),
$$
where $(x,y,z)$ is given by (\ref{determine}).

\noindent {\it Step 1.}\ \ By (\ref{423}), (\ref{51}) and (\ref{52}), we have
\begin{eqnarray}\label{OPP}
&&{\rm det}(\Gamma(x,y,z))\nonumber\\
&=&1-[f^{-1}(xy)]^2-[f^{-1}(xz)]^2-[f^{-1}(yz)]^2+2f^{-1}(xy)f^{-1}(xz)f^{-1}(yz)\nonumber\\
&=&1-\varsigma^2-\eta^2-\xi^2+2\varsigma\eta\xi\nonumber\\
&=&1-\cos^2\theta-\cos^2\nu-\cos^2\mu+2\cos\theta\cos\nu\cos\mu\nonumber\\
&=&\sin^2\mu\sin^2\nu-(\cos\theta-\cos\mu\cos\nu)^2.
\end{eqnarray}
Then,
\begin{eqnarray}\label{determinant2}
{\rm det}(\Gamma(x,y,z))>0\Leftrightarrow\left\{\begin{array}{ll}
\theta\in (|\mu-\nu|,\mu+\nu), &  \text{\ \ \ \ \ \ \ \ if }\mu+\nu\le\pi, \\
\theta\in (|\mu-\nu|,2\pi-(\mu+\nu)), &  \text{\ \ \ \ \ \ \ \ if }\mu+\nu>\pi.
\end{array}
\right.
\end{eqnarray}
Let $\mu,\nu\in(0,\pi)$. Define
$$
V_{\mu\nu}=\{\theta\in(0,\pi):\ {\rm condition}\ (\ref{determinant2})\ {\rm holds}\}.
$$
We will show that $J(\cdot,\nu,\mu)$ is an increasing function on $V_{\mu\nu}$.

By (\ref{423}), (\ref{WW5}) and (\ref{51}), we get
\begin{eqnarray}\label{determinant}
&&H(x,y,z)\nonumber\\
&=&\frac{1}{1-[f^{-1}(xy)]^2-[f^{-1}(xz)]^2-[f^{-1}(yz)]^2+2f^{-1}(xy)f^{-1}(xz)f^{-1}(yz)}\nonumber\\
&&\cdot\{x^2(1-[f^{-1}(yz)]^2)+y^2(1-[f^{-1}(xz)]^2)+z^2(1-[f^{-1}(xy)]^2)\nonumber\\
&&\ \ \ \ +2xy[f^{-1}(xz)f^{-1}(yz)-f^{-1}(xy)]+2xz[f^{-1}(xy)f^{-1}(yz)-f^{-1}(xz)]\nonumber\\
&&\ \ \ \ +2yz[f^{-1}(xy)f^{-1}(xz)-f^{-1}(yz)]\}\nonumber\\
&=&\frac{1}{1-\varsigma^2-\eta^2-\xi^2+2\varsigma\eta\xi}
\left\{\frac{f(\varsigma)f(\eta)(1-\xi^2)}{f(\xi)}+\frac{f(\varsigma)f(\xi)(1-\eta^2)}{f(\eta)}+\frac{f(\xi)f(\eta)(1-\varsigma^2)}{f(\varsigma)}\right.\nonumber\\
&&\left.\ \ \ \ +2f(\varsigma)(\eta\xi-\varsigma)+2f(\eta)(\varsigma\xi-\eta)+2f(\xi)(\varsigma\eta-\xi)\right\}.
\end{eqnarray}
For $\tau\in(-1,1)$, define
$$
\varsigma=\xi\eta+\tau\sqrt{1-\xi^2}\sqrt{1-\eta^2}.
$$
Then,
\begin{equation}\label{57}
\cos\theta=\cos\mu\cos\nu+ \tau\sin\mu\sin\nu,
\end{equation}
and hence
\begin{eqnarray}\label{RR1}
&&1-\varsigma^2-\eta^2-\xi^2+2\varsigma\eta\xi\nonumber\\
&=&1-\cos^2\theta-\cos^2\nu-\cos^2\mu+2\cos\theta\cos\nu\cos\mu\nonumber\\
&=&\sin^2\mu\sin^2\nu-(\cos\theta-\cos\mu\cos\nu)^2\nonumber\\
&=&(1-\tau^2)\sin^2\mu\sin^2\nu.
\end{eqnarray}

In this step, we assume without loss of generality that $\mu\ge\nu$. Then,
\begin{eqnarray*}
\left\{\begin{array}{ll}
\theta\in (\mu-\nu,\mu+\nu), &  \text{\ \ \ \ \ \ \ \ if }\mu+\nu\le\pi, \\
\theta\in (\mu-\nu,2\pi-(\mu+\nu)), &  \text{\ \ \ \ \ \ \ \ if }\mu+\nu>\pi.
\end{array}
\right.
\end{eqnarray*}
Denote
\begin{equation}\label{59}
a=\mu+\nu,\ \ \ \ b=\mu-\nu.
\end{equation}
By (\ref{422}), (\ref{52}) and (\ref{determinant})--(\ref{59}), we get
{\small\begin{eqnarray*}
&&J(\theta,\nu,\mu)\nonumber\\
&=&H(x,y,z)\nonumber\\
&=&\frac{1}{(1-\tau^2)\sin^2\mu\sin^2\nu}
\left\{f(\cos\theta)\left(\frac{\mu\sin\mu\sin\nu}{\nu}+\frac{\nu\sin\nu\sin\mu}{\mu}-{2\tau\sin\mu\sin\nu}\right)\right.\\
&&\ \ \ \ +\frac{\sin\mu\sin\nu(1-\cos^2\theta)}{\mu\nu f(\cos\theta)}\\
&&\ \ \ \ \left.+\frac{2\sin\nu\sin\mu(\tau\cos\mu\sin\nu-\sin\mu\cos\nu)}{\nu}
+\frac{2\sin\mu\sin\nu(\tau\sin\mu\cos\nu-\cos\mu\sin\nu)}{\mu}\right\}\\
&=&\frac{1}{(1-\tau^2)\mu\nu\sin\mu\sin\nu}
\left\{f(\cos\theta)\left(\mu^2+\nu^2-{2\tau\mu\nu}\right)+\frac{1-\cos^2\theta}{f(\cos\theta)}\right.\\
&&\ \ \ \ \left.+{2\mu(\tau\cos\mu\sin\nu-\sin\mu\cos\nu)}
+{2\nu(\tau\sin\mu\cos\nu-\cos\mu\sin\nu)}\right\}\\
&=&\frac{2[\cos(\mu-\nu)-\cos(\mu+\nu)] }{[(\mu+\nu)^2-(\mu-\nu)^2]\cdot[\cos\theta-\cos(\mu+\nu)]\cdot[\cos(\mu-\nu)-\cos\theta]}\\
&&\ \ \ \ \cdot\left\{\frac{\sin\theta}{\theta}\left[(\mu-\nu)^2+\frac{[(\mu+\nu)^2-(\mu-\nu)^2]\cdot[\cos(\mu-\nu)-\cos\theta]}{\cos(\mu-\nu)-\cos(\mu+\nu)}\right]
+\theta\sin\theta\right.\\
&&\ \ \ \ \left.+\frac{2[\cos\theta-\cos(\mu-\nu)]}{\cos(\mu-\nu)-\cos(\mu+\nu)}[(\mu+\nu)\sin(\mu+\nu)-(\mu-\nu)\sin(\mu-\nu)]\right.\\
&&\ \ \ \ \left.-2(\mu-\nu)\sin(\mu-\nu)\right\}\\
&=&\frac{2(\cos b-\cos a) }{(a^2-b^2)(\cos\theta-\cos a)(\cos b-\cos\theta)}\left\{\frac{\sin\theta}{\theta}\left[b^2+\frac{(a^2-b^2)(\cos b-\cos\theta)}{\cos b-\cos a}\right]
+\theta\sin\theta\right.\\
&&\ \ \ \ \left.+\frac{2(\cos\theta-\cos b)}{\cos b-\cos a}(a\sin a-b\sin b)-2 b\sin b\right\}\\
&=&\frac{2(\cos b-\cos a) }{(a^2-b^2)(\cos\theta-\cos a)(\cos b-\cos\theta)}\left\{\frac{\sin\theta}{\theta}\cdot\frac{a^2(\cos b-\cos\theta)+b^2(\cos \theta-\cos a)}{\cos b-\cos a}\right.\\
&&\ \ \ \ +\frac{\theta\sin\theta(\cos b-\cos\theta)}{\cos b-\cos a}+\frac{\theta\sin\theta(\cos\theta-\cos a)}{\cos b-\cos a}\\
&&\ \ \ \ \left.-\frac{2 a\sin a(\cos b-\cos\theta)}{\cos b-\cos a}-\frac{2 b\sin b(\cos\theta-\cos a)}{\cos b-\cos a}\right\}\\
&=&\frac{2 }{a^2-b^2}\left\{\frac{a^2\sin\theta+\theta^2\sin\theta-2a\theta\sin a }{\theta(\cos \theta-\cos a)}+\frac{b^2\sin\theta+\theta^2\sin\theta-2b\theta\sin b}{\theta(\cos b-\cos \theta)}\right\}\\
&=&\frac{2 }{a^2-b^2}\left\{\frac{(a-\theta)^2\sin\theta+2a\theta(\sin\theta-\sin a)}{\theta(\cos \theta-\cos a)}-\frac{(b-\theta)^2\sin\theta+2b\theta(\sin\theta-\sin b)}{\theta(\cos \theta-\cos b)}\right\}\\
&:=&\frac{2}{a^2-b^2}G(a,b,\theta).
\end{eqnarray*}}

We will show that $\partial_{\theta}G(a,b,\theta)>0$. For $u\in[0,\theta)\bigcup(\theta,2\pi-\theta)$, define
\begin{equation}\label{KKK1}
K(u,\theta)=\frac{(u-\theta)^2\sin \theta+2u \theta(\sin \theta-\sin u)}{\theta(\cos \theta-\cos u)}.
\end{equation}
Then,
$$
G(a,b,\theta)=K(a,\theta)-K(b,\theta),
$$
and
\begin{equation}\label{IN1}
\partial_{\theta}G(\theta,a,b)=\partial_{\theta}K(a,\theta)-\partial_{\theta}K(b,\theta).
\end{equation}
By (\ref{KKK1}), we get
{\small\begin{eqnarray}\label{IN2}
&&\partial_{\theta}K(u,\theta)\nonumber\\
&=&\frac{(u^2+\theta^2)\theta(1-\cos\theta\cos u)-(u^2-\theta^2)\sin\theta(\cos\theta-\cos u)-2 u\theta^2\sin u\sin\theta}{\theta^2(\cos\theta-\cos u)^2}\nonumber\\
&:=&P(u,\theta).
\end{eqnarray}}
By virtue of L'H\^opital's rule, we obtain that
{\small\begin{eqnarray*}
&&\lim_{u\rightarrow\theta}P(u,\theta)\nonumber\\
&=&\frac{1}{\theta^2}\lim_{u\rightarrow\theta}\frac{2u\theta(1-\cos\theta\cos u)-2u\sin\theta(\cos\theta-\cos u)-\sin u(u^2+\theta^2)(\sin\theta-\theta\cos\theta)-2u\theta^2\cos u\sin\theta}{2\sin u(\cos\theta-\cos u)}\\
&=&\frac{1}{\theta^2}\left\{-\theta+\lim_{u\rightarrow\theta}\frac{2u\theta(1-\cos\theta\cos u)-\sin u(u^2+\theta^2)(\sin\theta-\theta\cos\theta)-2u\theta^2\cos u\sin\theta}{2\sin u(\cos\theta-\cos u)}\right\}\\
&=&\frac{1}{\theta^2}\left\{-\theta+\lim_{u\rightarrow\theta}\frac{2u\theta(1-\cos\theta\cos u)-\sin u\sin\theta(u^2+\theta^2)}{2\sin \theta(\cos\theta-\cos u)}\right.\\
&&\ \ \ \ \ \ \ \ \left.+\lim_{u\rightarrow\theta}\frac{\theta\cos\theta\sin u(u^2+\theta^2)-2u\theta^2\cos u\sin\theta}{2\sin \theta(\cos\theta-\cos u)}\right\}\\
&=&\frac{1}{\theta^2}\left\{-\theta+\frac{\theta^3}{\sin^2\theta}\right\}\\
&=&\frac{\theta^2-\sin^2\theta}{\theta\sin^2\theta}\\
&:=&P(\theta).
\end{eqnarray*}}

We have
{\small\begin{eqnarray}\label{wan1}
&&P(u,\theta)>P(\theta)\nonumber\\
&\Leftrightarrow&\left[(u^2+\theta^2)\theta(1-\cos\theta\cos u)-(u^2-\theta^2)\sin\theta(\cos\theta-\cos u)-2 u\theta^2\sin u\sin\theta\right]\sin^2\theta\nonumber\\
&&>\theta(\theta^2-\sin^2\theta)(\cos\theta-\cos u)^2\nonumber\\
&\Leftrightarrow&\left\{(u^2+\theta^2)\theta\left[1-\cos(u-\theta)\right ]-(u^2-\theta^2)\sin\theta(\cos\theta-\cos u)+(u-\theta)^2\theta\sin u\sin\theta\right\}\sin^2\theta\nonumber\\
&&>\theta(\theta^2-\sin^2\theta)(\cos\theta-\cos u)^2\nonumber\\
&\Leftrightarrow&\left[2(u^2+\theta^2)\theta\sin^2\left(\frac{u-\theta}{2}\right)
-2(u^2-\theta^2)\sin\theta\sin\left(\frac{u-\theta}{2}\right)\sin\left(\frac{u+\theta}{2}\right)\right.\nonumber\\
&&\ \ \ \ \left.+(u-\theta)^2\theta\left(\sin^2\left(\frac{u+\theta}{2}\right)-\sin^2\left(\frac{u-\theta}{2}\right)\right)\right]\sin^2\theta\nonumber\\
&&>4\theta(\theta^2-\sin^2\theta)\sin^2\left(\frac{u-\theta}{2}\right)\sin^2\left(\frac{u+\theta}{2}\right)\nonumber\\
&\Leftrightarrow&\left[(u+\theta)^2\theta\sin^2\left(\frac{u-\theta}{2}\right)
-2(u^2-\theta^2)\sin\theta\sin\left(\frac{u-\theta}{2}\right)\sin\left(\frac{u+\theta}{2}\right)\right.\nonumber\\
&&\ \ \ \ \left.+(u-\theta)^2\theta\sin^2\left(\frac{u+\theta}{2}\right)\right]\sin^2\theta\nonumber\\
&&>4\theta(\theta^2-\sin^2\theta)\sin^2\left(\frac{u-\theta}{2}\right)\sin^2\left(\frac{u+\theta}{2}\right)\nonumber\\
&\Leftrightarrow&\left[\frac{(u+\theta)^2\theta}{\sin^2\left(\frac{u+\theta}{2}\right)}
-\frac{2(u^2-\theta^2)\sin\theta}{\sin\left(\frac{u-\theta}{2}\right)\sin\left(\frac{u+\theta}{2}\right)}
+\frac{(u-\theta)^2\theta}{\sin^2\left(\frac{u-\theta}{2}\right)}\right]\sin^2\theta>4\theta(\theta^2-\sin^2\theta)\nonumber\\
&\Leftrightarrow&\left[\frac{u+\theta}{\sin\left(\frac{u+\theta}{2}\right)}
-\frac{u-\theta}{\sin\left(\frac{u-\theta}{2}\right)}\right]^2\theta\sin^2\theta
+\frac{2(u^2-\theta^2)(\theta-\sin\theta)\sin^2\theta}{\sin\left(\frac{u-\theta}{2}\right)\sin\left(\frac{u+\theta}{2}\right)}>4\theta(\theta^2-\sin^2\theta).
\end{eqnarray}}

\noindent Case 1: Suppose that $u\in(\theta,2\pi-\theta)$.

Define
$$
g(w)=\frac{w}{\sin w},\ \ \ \ w\in(0,\pi).
$$
Then,
$$
g'(w)=\frac{\sin w- w\cos w}{\sin^2 w}>0,
$$
which implies that $g(w)$ is an increasing function for $w\in(0,\pi)$. Hence
$$
\frac{u^2-\theta^2}{\sin\left(\frac{u-\theta}{2}\right)\sin\left(\frac{u+\theta}{2}\right)}=
\frac{u+\theta}{\sin\left(\frac{u+\theta}{2}\right)}\cdot\frac{u-\theta}{\sin\left(\frac{u-\theta}{2}\right)}=4g\left(\frac{u+\theta}{2}\right)g
\left(\frac{u-\theta}{2}\right)
$$
is an increasing function for $u\in(\theta,2\pi-\theta)$.

We have
$$
g^{''}(w)=\frac{w(1+\cos^2w)- 2\sin w\cos w}{\sin^3 w}>0,
$$
which implies that $g'(w)$ is an increasing function for $w\in(0,\pi)$.
Hence
$$
\partial_u\left(\frac{\frac{u+\theta}{2}}{\sin\left(\frac{u+\theta}{2}\right)}\right)
-\partial_u\left(\frac{\frac{u-\theta}{2}}{\sin\left(\frac{u-\theta}{2}\right)}\right)=\frac{1}{2}\left[g'\left(\frac{u+\theta}{2}\right)-
g'\left(\frac{u-\theta}{2}\right)\right]>0,
$$
which implies that
$$
\frac{\frac{u+\theta}{2}}{\sin\left(\frac{u+\theta}{2}\right)}
-\frac{\frac{u-\theta}{2}}{\sin\left(\frac{u-\theta}{2}\right)}
$$
is an increasing function for $u\in(\theta,2\pi-\theta)$. Thus,
$$
\left[\frac{u+\theta}{\sin\left(\frac{u+\theta}{2}\right)}
-\frac{u-\theta}{\sin\left(\frac{u-\theta}{2}\right)}\right]^2\theta\sin^2\theta
+\frac{2(u^2-\theta^2)(\theta-\sin\theta)\sin^2\theta}{\sin\left(\frac{u-\theta}{2}\right)\sin\left(\frac{u+\theta}{2}\right)}
$$
is an increasing function for $u\in(\theta,2\pi-\theta)$.

\noindent Case 2:  Suppose that $u\in[0,\theta)$.

We have
$$
\frac{\frac{u+\theta}{2}}{\sin\left(\frac{u+\theta}{2}\right)}
-\frac{\frac{u-\theta}{2}}{\sin\left(\frac{u-\theta}{2}\right)}=
\frac{\frac{\theta+ u}{2}}{\sin\left(\frac{\theta+ u}{2}\right)}
-\frac{\frac{\theta -u}{2}}{\sin\left(\frac{\theta -u}{2}\right)}=g\left(\frac{u+\theta}{2}\right)-g
\left(\frac{\theta -u}{2}\right),
$$
which implies that
$$
\frac{\frac{u+\theta}{2}}{\sin\left(\frac{u+\theta}{2}\right)}
-\frac{\frac{u-\theta}{2}}{\sin\left(\frac{u-\theta}{2}\right)}
$$
is an increasing function for $u\in[0,\theta)$.

We have
\begin{eqnarray*}
(\ln g)^{''}(w)&=&\frac{g^{''}(w)g(w)-(g'(w))^2}{g^2(w)}\\
&=&\frac{w^2-\sin^2w}{g^2(w)\sin^4(w)}\\
&>&0,
\end{eqnarray*}
and
\begin{eqnarray*}
&&\partial_u\left(\frac{u^2-\theta^2}{\sin\left(\frac{u-\theta}{2}\right)\sin\left(\frac{u+\theta}{2}\right)}\right)>0\\
&\Leftrightarrow&
\partial_u\left(\frac{\frac{\theta+ u}{2}}{\sin\left(\frac{\theta+ u}{2}\right)}\cdot\frac{\frac{\theta- u}{2}}{\sin\left(\frac{\theta- u}{2}\right)}\right)>0\\
&\Leftrightarrow&g'\left(\frac{\theta+ u}{2}\right)g
\left(\frac{\theta -u}{2}\right)-g'
\left(\frac{\theta -u}{2}\right)g\left(\frac{\theta+ u}{2}\right)>0\\
&\Leftrightarrow&\frac{g'\left(\frac{\theta+ u}{2}\right)}{g\left(\frac{\theta+ u}{2}\right)}>
\frac{g'\left(\frac{\theta -u}{2}\right)}{g\left(\frac{\theta -u}{2}\right)}\\
&\Leftrightarrow&(\ln g)'\left(\frac{\theta +u}{2}\right)>(\ln g)'\left(\frac{\theta -u}{2}\right).
\end{eqnarray*}
Then,
$$
\frac{u^2-\theta^2}{\sin\left(\frac{u-\theta}{2}\right)\sin\left(\frac{u+\theta}{2}\right)}
$$
is an increasing function for $u\in[0,\theta)$. Thus,
$$
\left[\frac{u+\theta}{\sin\left(\frac{u+\theta}{2}\right)}
-\frac{u-\theta}{\sin\left(\frac{u-\theta}{2}\right)}\right]^2\theta\sin^2\theta
+\frac{2(u^2-\theta^2)(\theta-\sin\theta)\sin^2\theta}{\sin\left(\frac{u-\theta}{2}\right)\sin\left(\frac{u+\theta}{2}\right)}
$$
is an increasing function for $u\in[0,\theta)$.

By virtue of L'H\^opital's rule, we get
\begin{eqnarray}\label{LH}
&&\lim_{u\rightarrow\theta}\left\{\left[\frac{u+\theta}{\sin\left(\frac{u+\theta}{2}\right)}
-\frac{u-\theta}{\sin\left(\frac{u-\theta}{2}\right)}\right]^2\theta\sin^2\theta
+\frac{2(u^2-\theta^2)(\theta-\sin\theta)\sin^2\theta}{\sin\left(\frac{u-\theta}{2}\right)\sin\left(\frac{u+\theta}{2}\right)}\right\}\nonumber\\
&=&\left(\frac{2\theta}{\sin\theta}-2\right)^2\theta\sin^2\theta+8\theta\sin\theta(\theta-\sin\theta)\nonumber\\
&=&4\theta(\theta^2-\sin^2\theta).
\end{eqnarray}
Then, we obtain by (\ref{wan1}) and (\ref{LH}) that $P(u,\theta)>P(\theta)$ if $u\in(\theta,2\pi-\theta)$ and $P(u,\theta)<P(\theta)$ if $u\in[0,\theta)$. Thus, by (\ref{IN1}) and (\ref{IN2}), we get
\begin{equation}\label{rerr}
\partial_{\theta}G(\theta,a,b)=P(a,\theta)-P(b,\theta)=[P(a,\theta)-P(\theta)]+ [P(\theta)-P(b,\theta)]>0.
\end{equation}

\noindent {\it Step 2.}\ \ Let $x,y\in(0,\infty)$ satisfying $xy<1$. We will show that either $U_{xy}=\emptyset$ or there exist $z_1<z_2$ such that $U_{xy}=(z_1,z_2)$. By symmetry, we assume without loss of generality that $x\le y$ in this step. Then, $\mu\le\nu$.

By (\ref{OPP}), we get
$$
{\rm det}(\Gamma(x,y,z))=0\Leftrightarrow\mu=2\pi-(\nu+\theta),\ \theta-\nu,\ {\rm or}\ \nu-\theta.
$$
By (\ref{51}) and (\ref{52}), we have
\begin{equation}\label{nu1}
\frac{\sin\mu}{\mu}=\frac{y}{x}\frac{\sin\nu}{\nu}.
\end{equation}
If $\mu=2\pi-(\nu+\theta)$, by (\ref{nu1}), we get
\begin{equation}\label{dsd11}
\frac{\sin(2\pi-(\nu+\theta))}{2\pi-(\nu+\theta)}=\frac{y}{x}\frac{\sin\nu}{\nu}.
\end{equation}
Since $\frac{\sin w}{w}$ is a decreasing function for $w\in(0,\pi)$, there exists at most one $\nu$ satisfying (\ref{dsd11}).

If $\mu=\theta-\nu$, then $\nu<\theta$. By (\ref{nu1}), we get
\begin{equation}\label{dsd223}
\frac{\sin(\theta-\nu)}{\theta-\nu}=\frac{y}{x}\frac{\sin\nu}{\nu}.
\end{equation}
Since $\frac{\sin w}{w}$ is a decreasing function for $w\in(0,\pi)$, there exists at most one $\nu$ satisfying (\ref{dsd223}). Moreover, there exists a $\nu$ satisfying (\ref{dsd223})
only if
\begin{equation}\label{nec}
\frac{y}{x}\frac{\sin\theta}{\theta}<1.
\end{equation}

If $\mu=\nu-\theta$, then $\nu> \theta$. By (\ref{nu1}), we get
\begin{equation}\label{WWW}
\frac{\sin(\nu-\theta)}{\nu-\theta}=\frac{y}{x}\frac{\sin\nu}{\nu}.
\end{equation}
Define
$$
W(w)=\frac{w\sin(w-\theta)}{(w-\theta)\sin w},\ \ \ \ w\in[\theta,\pi).
$$
We have
$$
W'(w)=\frac{(\sin\theta)w(w-\theta)-\theta\sin w\sin(w-\theta)}{[(w-\theta)\sin w]^2}.
$$
Then,
$$
W'(w)>0\Leftrightarrow\frac{\sin(w-\theta)}{(w-\theta)}\cdot\frac{\sin w}{w}<\frac{\sin\theta}{\theta}.
$$
Since $\frac{\sin w}{w}$ is a decreasing function for $w\in(0,\pi)$, we get $W'(w)>0$ for $w\in(\theta,\pi)$. Thus, there exists at most one $\nu$ satisfying (\ref{WWW}). Moreover, since $W$ is an increasing function on $[\theta,\pi)$ and
$$
W(\theta)=\frac{\theta}{\sin\theta},
$$
there exists a $\nu$ satisfying (\ref{WWW})
only if
\begin{equation}\label{nec11}
\frac{y}{x}\frac{\sin\theta}{\theta}>1.
\end{equation}

By comparing conditions (\ref{nec}) and  (\ref{nec11}), we conclude that it is impossible that both $\mu=\theta-\nu$ and $\mu=\nu-\theta$ have solutions.
Hence, ${\rm det}(\Gamma(x,y,z))=0$ has at most two solutions. Note that
\begin{eqnarray*}
&&{\rm det}(\Gamma(x,y,0)):=-1-\varsigma^2+2\varsigma<0,\nonumber\\
&&{\rm det}\left(\Gamma\left(x,y,\frac{1}{y}\right)\right):=-\varsigma^2-\left[f^{-1}\left(\frac{x}{y}\right)\right]^2-2\varsigma f^{-1}\left(\frac{x}{y}\right)<0.
\end{eqnarray*}
If ${\rm det}(\Gamma(x,y,z))=0$ has at most one solution, then $U_{xy}=\emptyset$.
If ${\rm det}(\Gamma(x,y,z))=0$ has two solutions $z_1<z_2$, then we have that
$$
{\rm det}(\Gamma(x,y,z))<0\ \ \ \ {\rm for}\ z\in (0,z_1)\bigcup\left(z_2,\frac{1}{y}\right),
$$
and
$$
{\rm det}(\Gamma(x,y,z))>0,\ \ \forall z\in (z_1,z_2)\ \ {\rm or}\ \ {\rm det}(\Gamma(x,y,z))<0,\ \ \forall z\in (z_1,z_2).
$$
Hence ${\rm det}(\Gamma(x,y,z))>0$ has no solution or ${\rm det}(\Gamma(x,y,z))>0$ for any $z\in (z_1,z_2)$. Thus, $U_{xy}=\emptyset$ or $U_{xy}=(z_1,z_2)$.

\noindent {\it Step 3.}\ \ Define
$$
V=\{(\theta,\nu,\mu)\in(0,\pi)^3:\ {\rm condition}\ (\ref{determinant2})\ {\rm holds}\}.
$$
Then, $V$ is an open subset of $(0,\pi)^3$. Let $(\theta_0,\nu_0,\mu_0)\in V$. Then, there exists $\varepsilon>0$ such that
$$
V_0:=[\theta_0-\varepsilon,\theta_0+\varepsilon]\times [\nu_0-\varepsilon,\nu_0+\varepsilon]\times [\mu_0-\varepsilon,\mu_0+\varepsilon]\subset V.
$$
By Step 1 and symmetry of the function $J$, we conclude that for any $(\theta,\nu,\mu), (\theta',\nu',\mu')\in V_0$ satisfying $\theta\le \theta'$, $\nu\le \nu'$ and $\mu\le \mu'$, we have
$$
J(\theta,\nu,\mu)\le J(\theta',\nu',\mu'),
$$
and the equality holds if and only if $(\theta,\nu,\mu)=(\theta',\nu',\mu')$.

Define
$$
U=\{(x,y,z)\in(0,\infty)^3: xy,xz,yz<1,\ {\rm det}(\Gamma(x,y,z))>0\}.
$$
Let $(x_0,y_0,z_0)\in U$. Since $f$ is a homeomorphic map between $U$ and $V$, there exists $\epsilon>0$ such that for any $(x,y,z),(x',y',z')\in  [x_0-\epsilon,x_0+\epsilon]\times [y_0-\epsilon,y_0+\epsilon]\times [z_0-\epsilon,z_0+\epsilon]$ with $x\le x'$, $y\le y'$ and $z\le z'$, we have
$$
H(x,y,z)\ge H(x',y',z'),
$$
and the equality holds if and only if $(x,y,z)=(x',y',z')$.

Let $w_1,w_2\in (0,\infty)$ satisfying $w_1w_2<1$. Suppose that $w_3,w_4\in U_{w_1w_2}$ satisfying $w_3<w_4$. By Step 2, we know that there exist $z_1<z_2$ such that $U_{w_1w_2}=(z_1,z_2)$. By the above analysis, we know that for any $z\in (z_1,z_2)$, there exists $\epsilon_z>0$ such that for any $p,q\in [z-\epsilon_z,z+\epsilon_z]$ with $p<q$, we have
$$
H(w_1,w_2,p)>H(w_1,w_2,q).
$$
Therefore, the proof is complete by virtue of the Heine-Borel theorem.\hfill\fbox

\begin{rem}
To obtain the monotonicity of the function $H$, we prove in Step 1 the monotonicity of the function $G(\cdot,a,b)$  by establishing (\ref{rerr}). We believe that the following stronger result holds.
\begin{equation}\label{Appen1}
\partial_uP(u,\theta)>0,\ \ \ \ {\rm for\ any}\ u\in[0,\theta)\bigcup(\theta,2\pi-\theta)\ {\rm and}\ \theta\in(0,\pi).
\end{equation}
In the Appendix, we will show that (\ref{Appen1}) is equivalent to the following new inequality.
\begin{eqnarray}\label{Appen2}
&&(\cos\theta-\cos u)^2\left[2u\theta(1-\cos\theta\cos u)-2 u\sin\theta(\cos\theta-\cos u)+\sin u\sin\theta(u^2-3\theta^2)\right]\nonumber\\
&&\ \ +(\cos\theta-\cos u)\left[2u\theta^2\sin\theta(\sin^2u+1-\cos\theta\cos u)-(u^2+\theta^2)\theta\sin u(\sin^2\theta+1-\cos\theta\cos u)\right]\nonumber\\
&&>0,\ \ \ \ \ \ \ \ {\rm for\ any}\ u\in[0,\theta)\bigcup(\theta,2\pi-\theta)\ {\rm and}\ \theta\in(0,\pi).
\end{eqnarray}
So far we have not been able to prove (\ref{Appen2}), however, it has been verified by computer.
\end{rem}

\section{Some remarks}\setcounter{equation}{0}

In this section, we make some remarks. We hope they are helpful for considering the higher dimensional GRVM and SMW conjectures.

\subsection{Non-concavity of $F(\Lambda)$}
If $F$ is a concave function on $\bar{\cal S}$, then the sufficiency of Theorem \ref{main} can be established directly by Lemma \ref{Gauss2} and Theorem \ref{Formula}. Unfortunately, we find that $F$ is not a concave function.
Consider
\begin{eqnarray*}
\Lambda=\left(\begin{array}{cccc}
1&0.93&0.91&0.90\\
0.93&1&0.75&0.77\\
0.91&0.75&1&0.75\\
0.90&0.77&0.75&1
\end{array}\right),
\end{eqnarray*}
and ${\bar\Lambda}\in{\cal S}$ such that
$$({\bar\Lambda})_{kl}=\frac{\sum_{s<t}\Lambda_{st}}{6},\ \ k<l.
$$
By (\ref{June1}) and the software R, we get
$$
F(\Lambda)-F({\bar\Lambda})=0.0003994782,
$$
which implies that $F$ is non-concave.

We would like to take this opportunity to point out that numerical calculation via computer has played an important role in discovering some results of this paper.
Besides using Mathematica to quickly verify (\ref{equal2}) and using the software R to check non-concavity of $F$, we have used computer to check the monotonicity of $H$. It had greatly improved our confidence before we were able to give the proof of Theorem \ref{maint} at last.

We now give an affirmative result on the concavity of the function $F(\cdot)$ around some particular $\Lambda$. Let $X\in{\cal G}$. We denote by $S$ the  tetrahedron with vertexes $\{X_1,X_2,X_3,X_4\}$. If none of its dihedral angles is obtuse,
then $S$ is said to be nonobtuse. For $\Lambda\in {\cal S}$, we denote the Hessian matrix of $-F(\Lambda)$ by
$$
{\cal H}(\Lambda)=-\left(\frac{\partial^2 F}{\partial{\Lambda}_{pq}\partial{\Lambda}_{kl}}\right),\ \ \ \ p<q, k<l.
$$

\begin{pro} Suppose that $S$ is nonobtuse. Then, all diagonal elements of ${\cal H}(\Lambda^X)$ are positive and ${\rm det}({\cal H}(\Lambda^X))>0$.
\end{pro}
\noindent {\bf Proof.}\ \ Let $X\in{\cal G}$. To simplify notation, we denote $\Lambda^X$ by $\Lambda$. First, we show that
\begin{equation}\label{lam0}
-\frac{\partial^2 F}{\partial{\Lambda}_{kl}\partial{\Lambda}_{kl}}>0,\ \ \ \ \forall k<l.
\end{equation}

Denote by $\Upsilon^{X}_2$ the vertex Gramian
of $S$ associated with vertex $X_2$, which is equal to the covariance matrix of the Gaussian random vector $\left(\frac{X_1-X_2}{\sqrt{2}}, \frac{X_3-X_2}{\sqrt{2}},\frac{X_4-X_2}{\sqrt{2}}\right)$ as defined in (\ref{mmmm}). Since $S$ is nonobtuse, each facet of $S$ is nonobtuse (cf. \cite{F} and \cite[Proposition 2.8]{BC}). Hence the three facets meeting at vertex $X_2$
are all nonobtuse, which implies that (cf. (\ref{mmmm}) of the Appendix)
\begin{equation}\label{qiguai0}
r_1\ge 0,\ \ r_2\ge 0,\ \ r_3\ge 0.
\end{equation}
By \cite[Theorem 2.4]{BC}, $\Upsilon^{X}_2$ is a weakly diagonally dominant
Stieltjes matrix. Together with (\ref{zxcv}) of the Appendix, this implies that
\begin{equation}\label{qiguai}
{\tilde \Lambda}_{kl}\le 0,\ \ \ \ \forall k<l.
\end{equation}
Then, we obtain by (\ref{61}), (\ref{qiguai0}) and (\ref{qiguai}) that
$$
-\frac{\partial^2 F}{\partial{\Lambda}_{12}\partial{\Lambda}_{12}}>0.
$$
Similarly, we can show that (\ref{lam0}) holds for all $k<l$.

We now show that ${\rm det}({\cal H}(\Lambda))>0$. Denote by $\Phi(\Lambda)$ the $6\times 6$ diagonal matrix with
\begin{eqnarray*}
\Phi(\Lambda)_{11}&=&\frac{1}{8\sqrt{\pi^3(\Lambda')^3_{12}}}
\arccos\left(\frac{(\tilde{\Lambda})_{12}}{\sqrt{({\Lambda}')_{12}{\tilde A}^2+(\tilde{\Lambda})^2_{12}}}\right),\\
\Phi(\Lambda)_{22}&=&\frac{1}{8\sqrt{\pi^3(\Lambda')^3_{13}}}
\arccos\left(\frac{(\tilde{\Lambda})_{13}}{\sqrt{({\Lambda}')_{13}{\tilde A}^2+(\tilde{\Lambda})^2_{13}}}\right),\\
\cdots\cdots&&\\
\Phi(\Lambda)_{66}&=&\frac{1}{8\sqrt{\pi^3(\Lambda')^3_{34}}}
\arccos\left(\frac{(\tilde{\Lambda})_{34}}{\sqrt{({\Lambda}')_{34}{\tilde A}^2+(\tilde{\Lambda})^2_{34}}}\right).
\end{eqnarray*}
Define
\begin{equation}\label{phi}
\Psi(\Lambda)={\cal H}(\Lambda)-\Phi(\Lambda).
\end{equation}
Similar to (\ref{equal2}), we can show that
\begin{eqnarray}\label{dett}
\sum_{k<l}({\Lambda}')_{kl}\frac{\partial^2 F}{\partial \Lambda_{pq}\partial \Lambda_{kl}}=\frac{1}{2}\frac{\partial F}{\partial \Lambda_{pq}},
\ \ \ \ \forall p<q.
\end{eqnarray}
Define
$$v=(({\Lambda}')_{12},({\Lambda}')_{13},\dots,({\Lambda}')_{34}),$$
and
$$u=\left(\frac{\partial F}{\partial \Lambda_{12}},\frac{\partial F}{\partial \Lambda_{13}},\dots,\frac{\partial F}{\partial \Lambda_{34}}\right).
$$ Then, (\ref{dett}) can be written as
\begin{eqnarray}\label{psi}
-{\cal H}(\Lambda)v^T=\frac{1}{2}u^T.
\end{eqnarray}
By Theorem \ref{Gauss33}, Theorem \ref{Gauss44}, (\ref{phi}) and  (\ref{psi}), we get
$$
\Psi(\Lambda) v^T=0,
$$
which implies that ${\rm det}(\Psi(\Lambda))=0$. Therefore, we obtain by (\ref{lam0}) and (\ref{phi}) that
\begin{eqnarray*}
{\rm det}({\cal H}(\Lambda))&=&{\rm det}(\Phi(\Lambda))+{\rm det}(\Psi(\Lambda))+{\rm det}(\Phi(\Lambda))\cdot {\rm Trace}[\Phi(\Lambda)^{-1}\Psi(\Lambda)]\\
&=&{\rm det}(\Phi(\Lambda))\cdot\left\{1+{\rm Trace}[\Phi^{-1}(\Lambda)\Psi(\Lambda)]\right\}\\
&=&{\rm det}(\Phi(\Lambda))\cdot{\rm Trace}[\Phi^{-1}(\Lambda){\cal H}(\Lambda)]\\
&=&\left(\prod_{i=1}^6\Phi_{ii}(\Lambda)\right)\left(\sum_{i=1}^6\frac{{\cal H}_{ii}(\Lambda)}{\Phi_{ii}(\Lambda)}\right)\\
&>&0.
\end{eqnarray*}
\hfill\fbox

\subsection{Bounds for  ${\cal M}(X)$}

In this subsection, we derive the lower and upper bounds for
the function ${\cal M}(X)$.
\begin{pro} Let $X\in{\cal G}$ with covariance matrix $\Lambda$. Then,
\begin{eqnarray}\label{holds}
\frac{1}{4\sqrt{\pi}}\sum_{k<l}\sqrt{(\Lambda')_{kl}}\le E[\max(X_1,X_2,X_3,X_4)]\le\frac{1}{3\sqrt{\pi}}\sum_{k<l}\sqrt{(\Lambda')_{kl}}.
\end{eqnarray}
\end{pro}
\noindent {\bf Proof.}\ \  We define $X^I=\max(X_1,X_2,X_3,X_4)$, $X^{IV}=\min(X_1,X_2,X_3,X_4)$, and let
$X^{II}$ and $X^{III}$ be the second and third largest values of $X_1,X_2,X_3,X_4$, respectively. Note that
$$
E[X^{IV}]=-E[X^{I}],\ \ \ \ E[X^{III}]=-E[X^{II}].
$$

We have
{\small\begin{eqnarray}\label{a1}
&&E[X^{IV}]\nonumber\\
&=&E\left[X_1\wedge X_2\wedge X_3\wedge X_4\right]\nonumber\\
&=&E\left[(X_1-X_4)\wedge (X_2-X_4)\wedge (X_3-X_4)\wedge 0\right]\nonumber\\
&=&E\left[(X_1-X_4)\wedge (X_2-X_4)\wedge (X_3-X_4)\right]\nonumber\\
&&-E\left[\{(X_1-X_4)\wedge (X_2-X_4)\wedge (X_3-X_4)\}\vee 0\right]\nonumber\\
&=&E\left[X_1\wedge X_2\wedge X_3\right]-E\left[(X_1-X_4)^+\wedge (X_2-X_4)^+\wedge (X_3-X_4)^+\right]\nonumber\\
&=&E\left[X_1\wedge X_2\wedge X_3\right]-E\left[\{(X_1-X_4)\wedge (X_2-X_4)\wedge (X_3-X_4)\}\cdot 1_{\{X_4=\min_{1\leq i\leq 4}X_i\}}\right]\nonumber\\
&=&E\left[X_1\wedge X_2\wedge X_3\right]-E\left[\{X_1\wedge X_2\wedge X_3\}\cdot 1_{\{X_4=\min_{1\leq i\leq 4}X_i\}}\right]+E\left[X_4\cdot 1_{\{X_4=\min_{1\leq i\leq 4}X_i\}}\right]\nonumber\\
&=&E\left[X_1\wedge X_2\wedge X_3\right]-E\left[X^{III}\cdot 1_{\{X_4=\min_{1\leq i\leq 4}X_i\}}\right]+E\left[X^{IV}\cdot 1_{\{X_4=\min_{1\leq i\leq 4}X_i\}}\right].
\end{eqnarray}}
Similarly, we have
{\small\begin{eqnarray}\label{a2}
E(X^{IV})&=&E\left[X_2\wedge X_3\wedge X_4\right]-E\left[X^{III}\cdot 1_{\{X_1=\min_{1\leq i\leq 4}X_i\}}\right]+E\left[X^{IV}\cdot 1_{\{X_1=\min_{1\leq i\leq 4}X_i\}}\right],\ \ \ \ \ \ \
\end{eqnarray}
\begin{eqnarray}\label{a3}
E(X^{IV})&=&E\left[X_3\wedge X_4\wedge X_1\right]-E\left[X^{III}\cdot 1_{\{X_2=\min_{1\leq i\leq 4}X_i\}}\right]+E\left[X^{IV}\cdot 1_{\{X_2=\min_{1\leq i\leq 4}X_i\}}\right],\ \ \ \ \ \ \
\end{eqnarray}
\begin{eqnarray}\label{a4}
E(X^{IV})&=&E\left[X_4\wedge X_1\wedge X_2\right]-E\left[X^{III}\cdot 1_{\{X_3=\min_{1\leq i\leq 4}X_i\}}\right]+E\left[X^{IV}\cdot 1_{\{X_3=\min_{1\leq i\leq 4}X_i\}}\right].\ \ \ \ \ \ \
\end{eqnarray}}
By (\ref{ert0}), we get
\begin{eqnarray}\label{bb1}
E\left[X_1\wedge X_2\wedge X_3\right]=-\frac{1}{2\sqrt{\pi}}\left[\sqrt{(\Lambda')_{12}}+\sqrt{(\Lambda')_{13}}+\sqrt{(\Lambda')_{23}}\right],
\end{eqnarray}
\begin{eqnarray}\label{bb2}E\left[X_2\wedge X_3\wedge X_4\right]=-\frac{1}{2\sqrt{\pi}}\left[\sqrt{(\Lambda')_{23}}+\sqrt{(\Lambda')_{24}}+\sqrt{(\Lambda')_{34}}\right],
\end{eqnarray}
\begin{eqnarray}\label{bb3}E\left[X_3\wedge X_4\wedge X_1\right]=-\frac{1}{2\sqrt{\pi}}\left[\sqrt{(\Lambda')_{34}}+\sqrt{(\Lambda')_{13}}+\sqrt{(\Lambda')_{14}}\right],
\end{eqnarray}
\begin{eqnarray}\label{bb4}E\left[X_4\wedge X_1\wedge X_2\right]=-\frac{1}{2\sqrt{\pi}}\left[\sqrt{(\Lambda')_{14}}+\sqrt{(\Lambda')_{24}}+\sqrt{(\Lambda')_{12}}\right].
\end{eqnarray}

By summing up  (\ref{a1})--(\ref{a4}) and using (\ref{bb1})--(\ref{bb4}), we get
\begin{eqnarray*}
4E[X^{IV}]&=&E\left[X_1\wedge X_2\wedge X_3\right]+E\left[X_2\wedge X_3\wedge X_4\right]+E\left[X_3\wedge X_4\wedge X_1\right]+E\left[X_4\wedge X_1\wedge X_2\right]\nonumber\\
&&-E\left[X^{III}\right]+E\left[X^{IV}\right]\nonumber\\
&=&-\frac{1}{\sqrt{\pi}}\sum_{k<l}\sqrt{(\Lambda')_{kl}}-E\left[X^{III}\right]+E\left[X^{IV}\right].
\end{eqnarray*}
Then, we have that
\begin{eqnarray}\label{c1}
E\left[X^{III}\right]=3E\left[X^{I}\right]
-\frac{1}{\sqrt{\pi}}\sum_{k<l}\sqrt{(\Lambda')_{kl}},
\end{eqnarray}
and
\begin{eqnarray}\label{c2}
E\left[X^{II}\right]=-3E\left[X^{I}\right]
+\frac{1}{\sqrt{\pi}}\sum_{k<l}\sqrt{(\Lambda')_{kl}}.
\end{eqnarray}
Therefore, (\ref{holds}) holds by (\ref{c1}), (\ref{c2}) and the fact that $E\left[X^{III}\right]\le E\left[X^{II}\right]\le E\left[X^{I}\right]$.\hfill\fbox

\section{Appendix}\setcounter{equation}{0}

\subsection{Some auxiliary results of $\Lambda$}

Suppose that $X\in{\cal G}$ with covariance matrix $\Lambda\in {\cal S}$. Then, $\left(\frac{X_1-X_2}{\sqrt{2}}, \frac{X_3-X_2}{\sqrt{2}},\frac{X_4-X_2}{\sqrt{2}}\right)$ is a 3-dimensional centered
Gaussian random vector with covariance matrix
\begin{eqnarray}\label{mmmm}
\Upsilon^{X}_2&:=&\left(
\begin{array}{ccc}
1-{\Lambda}_{12}&\frac{1-{\Lambda}_{12}+{\Lambda}_{13}-{\Lambda}_{23}}{2} & \frac{1-{\Lambda}_{12}+{\Lambda}_{14}-{\Lambda}_{24}}{2}\\
\frac{1-{\Lambda}_{12}+{\Lambda}_{13}-{\Lambda}_{23}}{2} &1-{\Lambda}_{23}&\frac{1-{\Lambda}_{23}-{\Lambda}_{24}+{\Lambda}_{34}}{2}\\
\frac{1-{\Lambda}_{12}+{\Lambda}_{14}-{\Lambda}_{24}}{2}&\frac{1-{\Lambda}_{23}-{\Lambda}_{24}+{\Lambda}_{34}}{2}&1-{\Lambda}_{24}
\end{array}
\right)\nonumber\\
&:=&\left(
\begin{array}{ccc}
a_0&r_1 & r_2\\
r_1&b_0&r_3\\
r_2&r_3&c_0
\end{array}
\right).
\end{eqnarray}
We have
\begin{eqnarray*}
(\Upsilon^{X}_2)^{-1}&=&\frac{1}{\det(\Upsilon^{X}_2)}\left(
\begin{array}{ccc}
b_0c_0-r_3^2&r_2r_3-c_0r_1&r_1r_3-b_0r_2\\
r_2r_3-c_0r_1&a_0c_0-r_2^2&r_1r_2-a_0r_3\\
r_1r_3-b_0r_2&r_1r_2-a_0r_3&a_0b_0-r_1^2
\end{array}
\right)\nonumber\\
&:=&\frac{1}{\det(\Upsilon^{X}_2)}\left(
\begin{array}{ccc}
a&\rho_1&\rho_2\\
\rho_1&b&\rho_3\\
\rho_2&\rho_3&c
\end{array}
\right).
\end{eqnarray*}

By elementary calculation, we get
\begin{eqnarray}\label{zxcv}
&&\rho_1=\frac{1}{4}{\tilde \Lambda}_{23},\ \ \ \ \rho_2=\frac{1}{4}{\tilde \Lambda}_{24},\ \ \ \ \rho_3=\frac{1}{4}{\tilde \Lambda}_{12},\nonumber\\
&&a+\rho_1+\rho_2=-\frac{1}{4}{\tilde \Lambda}_{34},\nonumber\\
&&b+\rho_1+\rho_3=-\frac{1}{4}{\tilde \Lambda}_{13},\nonumber\\
&&c+\rho_2+\rho_3=-\frac{1}{4}{\tilde \Lambda}_{14}.
\end{eqnarray}
Then,
\begin{eqnarray}\label{ADD123}
({\Lambda}')_{12}{\tilde A}^2+(\tilde{\Lambda})^2_{12}
&=& 16\left[a_0\cdot{\rm det}(\Upsilon^{X}_2)+\rho^2_3\right]\nonumber\\
&=&16(a_0b_0-r_1^2)(a_0c_0-r_2^2)\nonumber\\
&=&\left\{4 ({\Lambda}')_{12}({\Lambda}')_{23}-[({\Lambda}')_{12}+({\Lambda}')_{23}-({\Lambda}')_{13}]^2 \right\}\nonumber\\
&&\ \ \ \ \cdot\left\{4 ({\Lambda}')_{12}({\Lambda}')_{24}-[({\Lambda}')_{12}+({\Lambda}')_{24}-({\Lambda}')_{14}]^2\right\}.
\end{eqnarray}
By
$$
(\tilde{\Lambda})_{12}=4(r_1r_2-a_0r_3),
$$
\begin{equation}\label{deg2}
4(\Lambda')_{12}(\Lambda')_{23}
-[(\Lambda')_{12}+(\Lambda')_{23}-(\Lambda')_{13}]^2=4(a_0c_0-r_2^2),
\end{equation}
and
\begin{equation}\label{deg3}
4(\Lambda')_{12}(\Lambda')_{24}
-[(\Lambda')_{12}+(\Lambda')_{24}-(\Lambda')_{14}]^2= 4(a_0b_0-r_1^2),
\end{equation}
we get
{\small\begin{eqnarray}\label{symm}
&&-2(\Lambda')_{12}+(\Lambda')_{13}+(\Lambda')_{14}
+(\Lambda')_{23}+(\Lambda')_{24}-2(\Lambda')_{34}\nonumber\\
&&\ \ +\frac{(\tilde{\Lambda})_{12}[(\Lambda')_{13}+(\Lambda')_{23}-(\Lambda')_{12}]}{4(\Lambda')_{12}(\Lambda')_{23}
-[(\Lambda')_{12}+(\Lambda')_{23}-(\Lambda')_{13}]^2}
+\frac{(\tilde{\Lambda})_{12}[(\Lambda')_{14}+(\Lambda')_{24}-(\Lambda')_{12}]}{4(\Lambda')_{12}(\Lambda')_{24}
-[(\Lambda')_{12}+(\Lambda')_{24}-(\Lambda')_{14}]^2}\nonumber\\
&=&4r_3-2r_1-2r_2+\frac{2(r_1r_2-a_0r_3)(c_0-r_2)}{a_0c_0-r_2^2}
+\frac{2(r_1r_2-a_0r_3)(b_0-r_1)}{a_0b_0-r_1^2}\nonumber\\
&=&\frac{2r_2(c_0r_1-r_2r_3+a_0r_3-r_1r_2-a_0c_0+r_2^2)}{a_0c_0-r_2^2}
+\frac{2r_1(b_0r_2-r_1r_3+a_0r_3-r_1r_2-a_0b_0+r_1^2)}{a_0b_0-r_1^2}\nonumber\\
&=&\frac{-2r_2(b+\rho_1+\rho_3)}{a_0c_0-r_2^2}
+\frac{-2r_1(c+\rho_2+\rho_3)}{a_0b_0-r_1^2}\nonumber\\
&=&\frac{(\tilde{\Lambda})_{13}[(\Lambda')_{12}+(\Lambda')_{23}-(\Lambda')_{13}]}{4(\Lambda')_{12}(\Lambda')_{23}
-[(\Lambda')_{12}+(\Lambda')_{23}-(\Lambda')_{13}]^2}
+\frac{(\tilde{\Lambda})_{14}[(\Lambda')_{12}+(\Lambda')_{24}-(\Lambda')_{14}]}{4(\Lambda')_{12}(\Lambda')_{24}
-[(\Lambda')_{12}+(\Lambda')_{24}-(\Lambda')_{14}]^2}.\ \ \ \ \ \
\end{eqnarray}}

\subsection{Calculation of a double integral}

In this subsection, we let $A>0$ and $$
\Psi=\left(\begin{array}{cc}
a_1&c_1\\
c_1&b_1
\end{array}\right)
$$
be a positive-definite matrix. Define
\begin{eqnarray*}
{\cal I}&=&\int_0^{\infty}\int_0^{\infty}e^{-\frac{a_1y^2+b_1z^2+2c_1yz}{2A^2}}dydz.
\end{eqnarray*}
\begin{lem}\label{lem71} We have
\begin{eqnarray}\label{nan2}
{\cal I}=\frac{A^2}{\sqrt{a_1b_1-c^2_1}}{\arccos\left(\frac{c_1}{\sqrt{a_1b_1}}\right)}.
\end{eqnarray}
\end{lem}
\noindent {\bf Proof.}\ \ The two eigenvalues of $\Psi$ are
$$
\lambda_{1}=\frac{a_1+b_1+\sqrt{(a_1-b_1)^2+4c_1^2}}{2},\ \ \ \ \lambda_{2}=\frac{a_1+b_1-\sqrt{(a_1-b_1)^2+4c_1^2}}{2}.
$$
Then, there exists an orthogonal matrix $O$ such that
\begin{eqnarray}\label{orthogonal-transformation}
O\left(\begin{array}{cc}
a_1&c_1\\
c_1&b_1
\end{array}\right)O^T=\left(\begin{array}{cc}
\lambda_{1}&0\\
0&\lambda_{2}
\end{array}\right).
\end{eqnarray}
Denote
$$
O=
\left(\begin{array}{cc}
\cos \theta&\sin\theta\\
-\sin\theta&\cos \theta
\end{array}\right),
$$
where $\theta\in (-\frac{\pi}{2},\frac{\pi}{2})$. One can check that if $c_1=0$, then $\theta=0$; if $c_1\neq 0$, then
$$
\tan\theta=\frac{\lambda_{1}-a_1}{c_1}=\frac{\sqrt{(a_1-b_1)^2+4c_1^2}-(a_1-b_1)}{2c_1},
$$
which implies  that if $c_1>0$, then $\theta\in (0,\frac{\pi}{2})$; if $c_1<0$, then
$\theta\in (-\frac{\pi}{2},0)$.

Define
\begin{eqnarray*}
Q=\left(
\begin{array}{cc}
\sqrt{\lambda_{2}}&0\\
0&\sqrt{\lambda_{1}}
\end{array}\right)O
=\left(\begin{array}{cc}
\sqrt{\lambda_2}\cos \theta&\sqrt{\lambda_2}\sin\theta\\
-\sqrt{\lambda_1}\sin\theta&\sqrt{\lambda_1}\cos \theta
\end{array}\right).
\end{eqnarray*}
By (\ref{orthogonal-transformation}), we get
\begin{eqnarray*}
Q\left(\begin{array}{cc}
a_1&c_1\\
c_1&b_1
\end{array}\right)Q^T&=&
\left(
\begin{array}{cc}
\sqrt{\lambda_{2}}&0\\
0&\sqrt{\lambda_{1}}
\end{array}\right)O\left(\begin{array}{cc}
a_1&c_1\\
c_1&b_1
\end{array}\right)O^T\left(
\begin{array}{cc}
\sqrt{\lambda_{2}}&0\\
0&\sqrt{\lambda_{1}}
\end{array}\right)\\
&=&\left(
\begin{array}{cc}
\sqrt{\lambda_{2}}&0\\
0&\sqrt{\lambda_{1}}
\end{array}\right)\left(\begin{array}{cc}
\lambda_{1}&0\\
0&\lambda_{2}
\end{array}\right)\left(
\begin{array}{cc}
\sqrt{\lambda_{2}}&0\\
0&\sqrt{\lambda_{1}}
\end{array}\right)\\
&=&\left(\begin{array}{cc}
\lambda_1\lambda_2&0\\
0&\lambda_1\lambda_2
\end{array}\right).
\end{eqnarray*}
Define
$
\left(\begin{array}{c}
y\\
z
\end{array}\right)=Q^T\left(\begin{array}{c}
u\\
v
\end{array}\right)
$. Then
\begin{eqnarray*}
{\cal I}=\sqrt{\lambda_1\lambda_2}\iint_{D_Q}e^{-\frac{\lambda_1\lambda_2(u^2+v^2)}{2A^2}}dudv,
\end{eqnarray*}
where
\begin{eqnarray*}
D_Q=\left\{(u,v)\in \mathbb{R}^2: 0\leq \sqrt{\lambda_2}(\cos\theta)u-\sqrt{\lambda_1}(\sin\theta)v,\ 0\le
\sqrt{\lambda_2}(\sin\theta)u+\sqrt{\lambda_1}(\cos\theta)v\right\}.
\end{eqnarray*}

If $c_1=0$, then $\theta=0$ and hence $
D_Q=\left\{(u,v)\in \mathbb{R}^2: u\geq 0, v\geq 0 \right\}.
$ We have
\begin{eqnarray*}
{\cal I}&=&\sqrt{\lambda_1\lambda_2}\int_0^{\frac{\pi}{2}}\int_0^{\infty}re^{-\frac{\lambda_1\lambda_2r^2}{2A^2}}drd\phi\nonumber\\
&=&\frac{\pi A^2}{2\sqrt{\lambda_1\lambda_2}}.
\end{eqnarray*}
If $c_1>0$, then $\theta\in (0,\frac{\pi}{2})$. We have
\begin{eqnarray*}
D_Q&=&\left\{(u,v)\in \mathbb{R}^2: -\sqrt{\frac{\lambda_2}{\lambda_1}}\cdot({\tan\theta})u\le v\leq \sqrt{\frac{\lambda_2}{\lambda_1}}\cdot\frac{u }{\tan\theta}\right\}.
\end{eqnarray*}
Then,
{\small\begin{eqnarray*}
{\cal I}&=&\sqrt{\lambda_1\lambda_2}\int_{\arctan\left(-\sqrt{\frac{\lambda_2}{\lambda_1}}{\tan\theta}\right)}^{\arctan\left(\sqrt{\frac{\lambda_2}{\lambda_1}}\frac{1}{\tan\theta}\right)}
\int_0^{\infty}re^{-\frac{\lambda_1\lambda_2r^2}{2A^2}}drd\phi\nonumber\\
&=&\frac{A^2}{\sqrt{\lambda_1\lambda_2}}{\arctan\left(\frac{\sqrt{\lambda_1\lambda_2}}{(\lambda_1-\lambda_2)\cos\theta\sin\theta}\right)}\\
&=&\frac{A^2}{\sqrt{a_1b_1-c^2_1}}{\arctan\left(\frac{\sqrt{a_1b_1-c_1^2}}{c_1}\right)}\\
&=&\frac{A^2}{\sqrt{a_1b_1-c^2_1}}{\arccos\left(\frac{c_1}{\sqrt{a_1b_1}}\right)}.
\end{eqnarray*}}
If  $c_1<0$, then $\theta\in (-\frac{\pi}{2},0)$. We have
\begin{eqnarray*}
D_Q&=&\left\{(u,v)\in \mathbb{R}^2:-\sqrt{\frac{\lambda_2}{\lambda_1}}\cdot({\tan\theta})u\le v,\ \sqrt{\frac{\lambda_2}{\lambda_1}}\cdot\frac{u }{\tan\theta}\le v\right\}.
\end{eqnarray*}
Then,
{\small\begin{eqnarray*}
{\cal I}&=&\sqrt{\lambda_1\lambda_2}\int^{\pi+\arctan\left(\sqrt{\frac{\lambda_2}{\lambda_1}}\frac{1}{\tan\theta}\right)}_{\arctan\left(-\sqrt{\frac{\lambda_2}{\lambda_1}}{\tan\theta}\right)}
\int_0^{\infty}re^{-\frac{\lambda_1\lambda_2r^2}{2A^2}}drd\phi\nonumber\\
&=&\frac{A^2}{\sqrt{\lambda_1\lambda_2}}\left\{\pi-{\arctan\left(\frac{-\sqrt{\lambda_1\lambda_2}}{(\lambda_1-\lambda_2)\cos\theta\sin\theta}\right)}\right\}\\
&=&\frac{A^2}{\sqrt{a_1b_1-c^2_1}}\left\{\pi-{\arctan\left(\frac{-\sqrt{a_1b_1-c_1^2}}{c_1}\right)}\right\}\\
&=&\frac{A^2}{\sqrt{a_1b_1-c^2_1}}{\arccos\left(\frac{c_1}{\sqrt{a_1b_1}}\right)}.
\end{eqnarray*}}
Therefore, (\ref{nan2}) holds for any case.\hfill\fbox

\subsection{Mathematica code for verifying (\ref{equal2})}

\begin{verbatim}
In[1]:= Expand[((1 - b)^2 - (1 - b)*(1 - a + 1 - c + 1 - d + 1 - f -
        2*(1 - e)) + (a - d)*(c - f))*(1 -
     a + (1 - d) - (1 -
       b))*(4*(1 - a)*(1 - e) - (1 - a - e + c)^2) + ((1 - c)^2 - (1 -
         c)*(1 - a + 1 - b + 1 - e + 1 - f - 2*(1 - d)) + (a - e)*(b -
         f))*(1 -
     a + (1 - e) - (1 - c))*(4*(1 - a)*(1 - d) - (1 - a - d + b)^2) -
  2*((1 - d)^2 - (1 - d)*(1 - a + 1 - b + 1 - e + 1 - f -
        2*(1 - c)) + (a - b)*(e - f))*(1 -
     b)*(4*(1 - a)*(1 - e) - (1 - a - e + c)^2) -
  2*((1 - e)^2 - (1 - e)*(1 - a + 1 - c + 1 - d + 1 - f -
        2*(1 - b)) + (a - c)*(d - f))*(1 -
     c)*(4*(1 - a)*(1 - d) - (1 - a - d + b)^2) -
  2*((1 - b)^2 - (1 - b)*(1 - a + 1 - c + 1 - d + 1 - f -
        2*(1 - e)) + (a - d)*(c - f))*(1 -
     d)*(4*(1 - a)*(1 - e) - (1 - a - e + c)^2) -
  2*((1 - c)^2 - (1 - c)*(1 - a + 1 - b + 1 - e + 1 - f -
        2*(1 - d)) + (a - e)*(b - f))*(1 -
     e)*(4*(1 - a)*(1 - d) - (1 - a - d + b)^2) -
  2*(1 - f)*(4*(1 - a)*(1 - d) - (1 - a - d + b)^2)*(4*(1 - a)*(1 -
        e) - (1 - a - e + c)^2)]
Out[1]= 0
\end{verbatim}

\subsection{Derivation of inequality (\ref{Appen2})}

Let $K(u,\theta)$ and $P(u,\theta)$ be defined as in (\ref{KKK1}) and (\ref{IN2}), respectively. We have
{\small\begin{eqnarray}\label{B1}
&&\partial_uP(u,\theta)\nonumber\\
&=&\frac{1}{\theta^2(\cos\theta-\cos u)^4}\left\{(\cos\theta-\cos u)^2[2u\theta(1-\cos\theta\cos u)-2u\sin\theta(\cos\theta-\cos u)\right.\nonumber\\
&&\ \ \ \ -\sin u(u^2+\theta^2)(\sin\theta-\theta\cos\theta)-2u\theta^2\cos u\sin\theta]\nonumber\\
&&\ \ \left.+2\sin u(\cos\theta-\cos u)[(u^2-\theta^2)\sin\theta(\cos\theta-\cos u)+2u\theta^2\sin u\sin\theta-(u^2+\theta^2)\theta(1-\cos\theta\cos u)]\right\}\nonumber\\
&=&\frac{1}{\theta^2(\cos\theta-\cos u)^4}\left\{(\cos\theta-\cos u)^2\left[ \left(2u\theta(1-\cos\theta\cos u)+2\sin u\sin\theta(u^2-\theta^2)\right)\right.\right.\nonumber\\
&&\ \ \ \ \ \ \ \ \left.-\left(2u\sin\theta(\cos\theta-\cos u)+\sin u(u^2+\theta^2)(\sin\theta-\theta\cos\theta)+2u\theta^2\cos u\sin\theta\right)\right]\nonumber\\
&&\left.\ \ +(\cos\theta-\cos u)\left[4u\theta^2\sin^2 u\sin\theta-2(u^2+\theta^2)\theta\sin u(1-\cos\theta\cos u)\right]\right\}.
\end{eqnarray}}
On the other hand, by elementary calculation, we get
{\small\begin{eqnarray*}
\partial_{u}K(u,\theta)=\frac{2u\theta(1-\cos\theta\cos u)-(u^2+\theta^2)\sin\theta\sin u+2(\cos\theta-\cos u)(u\sin\theta-\theta\sin u)}{\theta(\cos\theta-\cos u)^2},
\end{eqnarray*}}
and
{\small\begin{eqnarray}\label{B2}
&&\partial_{\theta}\partial_{u}K(u,\theta)\nonumber\\
&=&\frac{1}{\theta^2(\cos\theta-\cos u)^4}\left\{(\cos\theta-\cos u)^2\left[ \sin u(u^2+\theta^2)(\sin\theta-\theta\cos\theta)-2u\sin\theta(\cos\theta-\cos u)\right.\right.\nonumber\\
&&\ \ \ \ \ \ \ \ \left.-4\theta^2\sin\theta\sin u+2u\theta(\theta\cos u\sin\theta+1-\cos\theta\cos u)\right]\nonumber\\
&&\left.\ \ +(\cos\theta-\cos u)[4u\theta^2\sin\theta(1-\cos\theta\cos u)-2\theta\sin^2\theta\sin u(u^2+\theta^2)]\right\}.
\end{eqnarray}}
Then, we obtain by (\ref{IN2}), (\ref{B1}) and (\ref{B2}) that
{\small\begin{eqnarray*}
&&\partial_uP(u,\theta)\nonumber\\
&=&\partial_{u}\partial_{\theta}K(u,\theta)\nonumber\\
&=&\frac{\partial_uP(u,\theta)+\partial_{\theta}\partial_{u}K(u,\theta)}{2}\\
&=&\frac{1}{\theta^2(\cos\theta-\cos u)^4}\left\{(\cos\theta-\cos u)^2\left[ 2u\theta(1-\cos\theta\cos u)-2 u\sin\theta(\cos\theta-\cos u)\right.\right.\\
&&\ \ \ \ \ \ \ \ \ \ \ \ \ \ \ \ \ \ \ \ \ \ \ \ \ \ \ \ \ \ \ \ \left.+\sin u\sin\theta(u^2-3\theta^2)\right]\\
&&\left.\ \ +(\cos\theta-\cos u)[2u\theta^2\sin\theta(\sin^2u+1-\cos\theta\cos u)-(u^2+\theta^2)\theta\sin u(\sin^2\theta+1-\cos\theta\cos u)]\right\}.
\end{eqnarray*}}
\vskip -0.6cm
\noindent Therefore, $\partial_uP(u,\theta)>0$ for any $u\in[0,\theta)\bigcup(\theta,2\pi-\theta)$ and $\theta\in(0,\pi)$ is equivalent to (\ref{Appen2}) holds.

\bigskip

{ \noindent {\bf\large Acknowledgments}\ \ \ \ This work was supported by the Natural Sciences and Engineering Research Council of Canada, and the National Natural Science Foundation of China (Grant No. 11771309, No. 11871184 and No. 61973096).

\end{document}